\documentclass[reqno,oneside,12pt]{amsproc}
\usepackage[unicode]{hyperref} 
\usepackage{fancybox,ascmac,comment,xcolor}
\usepackage{amssymb}
\usepackage{amsfonts,mathrsfs}
\usepackage{amsmath}
\usepackage{amsthm}
\usepackage{mathtools}
\usepackage{enumerate}
\usepackage{mleftright}
\DeclareMathAlphabet{\mathpzc}{OT1}{pzc}{m}{it}
\mleftright
\usepackage[all]{xy}
\calclayout
\mathtoolsset{showonlyrefs=true}
\numberwithin{equation}{section}
\makeatletter
\newcommand{\shortmathcal}[1]{\@tfor\ch:=#1\do{
\expandafter\edef\csname c\ch\endcsname{\noexpand\mathcal{\ch}}
}}
\newcommand{\shortmathbb}[1]{\@tfor\ch:=#1\do{
\expandafter\edef\csname bb\ch\endcsname{\noexpand\mathbb{\ch}}
}}
\newcommand{\shortmathbf}[1]{\@tfor\ch:=#1\do{
\expandafter\edef\csname b\ch\endcsname{\noexpand\mathbf{\ch}}
}}
\newcommand{\shortboldsymbol}[1]{\@tfor\ch:=#1\do{
\expandafter\edef\csname bs\ch\endcsname{\noexpand\boldsymbol{\ch}}
}}
\newcommand{\shortmathfrak}[1]{\@tfor\ch:=#1\do{
\expandafter\edef\csname f\ch\endcsname{\noexpand\mathfrak{\ch}}}}
\newcommand{\shortmathscr}[1]{\@tfor\ch:=#1\do{
\expandafter\edef\csname s\ch\endcsname{\noexpand\mathscr{\ch}}}}
\newcommand{\shortmathrm}[1]{\@tfor\ch:=#1\do{
\expandafter\edef\csname r\ch\endcsname{\noexpand\mathrm{\ch}}
}}
\makeatother
\shortmathcal{ABCDEFGHIJKLMNOPQRSTUVWXYZ}
\shortmathbb{ABCDEFGHIJKLMNOPQRSTUVWXYZ}
\shortmathbf{ABCDEFGHIJKLMNOPQRSTUVWXYZabcdefghijklmnopqrstuvwxyz}
\shortboldsymbol{ABCDEFGHIJKLMNOPQRSTUVWXYZabcdefghijklmnopqrstuvwxyz}
\shortmathfrak{ABCDEFGHIJKLMNOPQRSTUVWXYZabcdefghjklmnopqrstuvwxyz}
\shortmathscr{ABCDEFGHIJKLMNOPQRSTUVWXYZ}
\shortmathrm{ABCDEFGHIJKLMNOPQRSTUVWXYZabcdefghijklmnopqrstuvwxyz}
\DeclareFontEncoding{OT2}{}{}
\DeclareFontSubstitution{OT2}{cmr}{m}{n}
\DeclareFontFamily{OT2}{cmr}{\hyphenchar\font45 }
\DeclareFontShape{OT2}{cmr}{m}{n}{
<5><6><7><8><9>gen*wncyr
<10><10.95><12><14.4><17.28><20.74><24.88>wncyr10}{}
\DeclareFontShape{OT2}{cmr}{b}{n}{
<5><6><7><8><9>gen*wncyb
<10><10.95><12><14.4><17.28><20.74><24.88>wncyb10}{}
\DeclareMathAlphabet{\mathcyr}{OT2}{cmr}{m}{n}
\DeclareMathAlphabet{\mathcyb}{OT2}{cmr}{b}{n}
\SetMathAlphabet{\mathcyr}{bold}{OT2}{cmr}{b}{n}
\newcommand{\red}[1]{\textcolor{red}{#1}}

\newcommand{\Hom}{\mathrm{Hom}}
\newcommand{\id}{\mathrm{id}}

\newcommand{\CAlg}{\mathsf{CAlg}}

\newcommand{\Lie}{\mathrm{Lie}}
\newcommand{\Spec}{\mathrm{Spec}}

\newcommand{\Ad}{\mathrm{Ad}}
\newcommand{\Aut}{\mathrm{Aut}}

\newcommand{\hotimes}{\mathbin{\widehat{\otimes}}}
\newcommand{\Ker}{\mathrm{Ker}}

\newcommand{\emp}{\varnothing}
\newcommand{\ep}{\varepsilon}
\newcommand{\jump}[1]{\ensuremath{[\![#1]\!]}}
\newcommand{\ncjump}[1]{\ensuremath{\langle\!\langle #1\rangle\!\rangle}}
\renewcommand{\mod}{\mathrm{mod}}
\DeclareMathOperator{\Span}{span}

\newcommand{\sh}{\mathbin{\mathcyr{sh}}}
\newcommand{\Sh}{\mathpzc{Sh}}

\newcommand{\BIMU}{\mathsf{BIMU}}

\newcommand{\LU}{\mathsf{LU}}
\newcommand{\mmu}{\mathsf{mu}}

\newcommand{\preari}{\mathsf{preari}}
\newcommand{\preira}{\mathsf{preira}}

\newcommand{\GARI}{\mathsf{GARI}}

\newcommand{\ARI}{\mathsf{ARI}}

\newcommand{\adari}{\mathsf{adari}}

\newcommand{\logmu}{\mathsf{logmu}}

\newcommand{\invmu}{\mathsf{invmu}}

\renewcommand{\Pi}{\mathsf{Pi}}

\newcommand{\fos}{\mathfrak{os}}

\newcommand{\pal}{\mathsf{pal}}

\newcommand{\as}{\mathsf{as}}
\newcommand{\al}{\mathsf{al}}
\newcommand{\fol}{\mathfrak{ol}}

\newcommand{\il}{\mathsf{il}}
\newcommand{\Esena}{\fE\text{-}\mathsf{sena}}
\newcommand{\senu}{\mathsf{senu}}

\newcommand{\Opusnu}{\fO\text{-}\mathsf{pusnu}}

\newcommand{\pusinu}{\mathsf{pusinu}}

\newcommand{\ful}[1]{{}_{#1}\lceil}
\newcommand{\fll}[1]{{}_{#1}\lfloor}
\newcommand{\fur}[1]{\rceil_{#1}}
\newcommand{\flr}[1]{\rfloor_{#1}}

\renewcommand{\neg}{\mathsf{neg}}
\newcommand{\swap}{\mathsf{swap}}
\newcommand{\anti}{\mathsf{anti}}
\newcommand{\pari}{\mathsf{pari}}
\newcommand{\push}{\mathsf{push}}
\newcommand{\pus}{\mathsf{pus}}

\newcommand{\pusnu}{\mathsf{pusnu}}

\newcommand{\mantar}{\mathsf{mantar}}

\newcommand{\der}{\mathsf{der}}

\newcommand{\mi}{\mathsf{mi}}

\newcommand{\weak}{\mathsf{weak}}

\newcommand{\DMR}{\mathsf{DMR}}
\newcommand{\KRV}{\mathsf{KRV}}

\newcommand{\dmr}{\mathfrak{dmr}}
\newcommand{\krv}{\mathfrak{krv}}

\newcommand{\KZ}{\mathrm{KZ}}
\newcommand{\hyp}{\mathchar`-}

\newcommand{\hcV}{\widehat{\cV}}

\newcommand{\rev}{\mathrm{rev}}
\newcommand{\inert}{\mathrm{inert}}

\newcommand{\ser}{\mathsf{ser}}
\newcommand{\com}{\mathsf{com}}
\newcommand{\Lau}{\mathsf{Lau}}
\newcommand{\ma}{\mathsf{ma}}

\newcommand{\EF}{\mathsf{EF}}
\newcommand{\tr}{\mathfrak{tr}}

\renewcommand{\div}{\mathsf{div}}

\newcommand{\user}{\bu\text{-}\ser}
\newcommand{\ucom}{\bu\text{-}\com}
\newcommand{\uLau}{\bu\text{-}\Lau}
\newcommand{\vcst}{\bv\text{-}\mathsf{cst}}

\newcommand{\hff}{\hat{\ff}}

\newcommand{\hcW}{\widehat{\cW}}
\newcommand{\SolKV}{\mathsf{SolKV}}
\newcommand{\WDMR}{\mathsf{WDMR}}
\newcommand{\wdmr}{\mathfrak{wdmr}}
\newcommand{\aut}{\mathrm{aut}}
\newcommand{\TS}{\mathbf{TS}}

\newcommand{\tsh}{\mathbin{\widetilde{\sh}}}

\newcommand{\hzeta}{\widehat{\zeta}}
\newcommand{\Red}{\mathrm{red}}
\newcommand{\HK}{\mathrm{HK}}
\newcommand{\Ra}{\mathrm{Ra}}
\newcommand{\Eul}[1]{#1\text{-}\mathsf{Eul}}

\newcommand{\baral}{\underline{\al}}
\newcommand{\baril}{\underline{\il}}
\newcommand{\barfol}{\underline{\fol}}
\newcommand{\bch}{\mathrm{bch}}

\newcommand{\SAut}{\mathrm{SAut}}
\newcommand{\sder}{\mathfrak{sder}}
\newcommand{\cent}{\mathrm{cent}}
\newcommand{\baras}{\underline{\as}}
\newcommand{\barfos}{\underline{\fos}}
\newcommand{\ola}[1]{\overleftarrow{#1}}

\newtheorem{theorem}{Theorem}[section]
\newtheorem{proposition}[theorem]{Proposition}
\newtheorem{lemma}[theorem]{Lemma}
\newtheorem{corollary}[theorem]{Corollary}
\newtheorem{conjecture}[theorem]{Conjecture}
\theoremstyle{definition}
\newtheorem{definition}[theorem]{Definition}
\theoremstyle{remark}
\newtheorem{remark}[theorem]{Remark}

\allowdisplaybreaks[0]
\title[KV relation for MZVs]{Kashiwara--Vergne relation for multiple zeta values}
\author{Hanamichi Kawamura, Takumi Maesaka}
\address[Hanamichi Kawamura]{Department of Mathematics, Faculty of Science Division I, Tokyo University of Science, 1-3 Kagurazaka, Shinjuku-ku, Tokyo, 162-8601, Japan}
\email{1125512@ed.tus.ac.jp}
\address[Takumi Maesaka]{Department of Mathematics, Kyushu University, Motooka 744, Nishi-ku, Fukuoka, 819-0395, Japan}
\email{nozaki.takumi.912@s.kyushu-u.ac.jp}
\begin{document}
\begin{abstract}
    We give a explicit formulation of the Kashiwara--Vergne structure as a family of linear relations among multiple zeta values. According to some classical conjectures, this family is a candidate to be exhaustive. Moreover, we show a relationship between the first Kashiwara--Vergne equation and multitangent functions.
\end{abstract}
\maketitle
\setcounter{tocdepth}{2}
\tableofcontents
\setcounter{section}{-1}
\section{Introduction}\label{sec:introduction}
\subsection{Motivation}
The two groups, the \emph{double shuffle group} $\DMR_{0}$ and \emph{Kashiwara--Vergne group} $\KRV_{2}$, in the de Rham setting, are conjectured to be isomorphic. Since the $\DMR_{0}$ was made as a kind of formalization of \emph{multiple zeta values} (\emph{MZVs}), it is natural to expect that also the Kashiwara--Vergne group has such an interpretation. In this paper, we show that the Kashiwara--Vergne group $\KRV_{2}$ can be viewed as an ``weakened version'' of $\DMR_{0}$, as well as their torsors handling linear relations among multiple zeta values without modulo $\pi^{2}$.

\subsection{Integrals and series}
Consider the iterated integral symbol defined as
\[I(0;a_{1},\ldots,a_{k};1)\coloneqq\int_{0<t_{1}<\cdots<t_{k}<1}\prod_{i=1}^{k}\frac{dt_{i}}{t_{i}-a_{i}},\]
for elements $a_{1},\ldots,a_{k}\in\{0,1\}$ satisfying $(a_{1},a_{k})=(1,0)$.
The following identity defines MZVs in two ways, namely as nested series and iterated integrals.

\begin{theorem}
    Let $r,k_{1},\ldots,k_{r}$ be positive integers and assume $k_{r}\ge 2$. Then the nested series
    \[\zeta(k_{1},\ldots,k_{r})\coloneqq\sum_{0<n_{1}<\cdots<n_{r}}\frac{1}{n_{1}^{k_{1}}\cdots n_{r}^{k_{r}}}\]
    converges absolutely and coincides with the following iterated integral:
    \[\zeta(k_{1},\ldots,k_{r}) = (-1)^{r}I(0;1,\underbrace{0,\ldots,0}_{k_{1}-1},\ldots,1,\underbrace{0,\ldots,0}_{k_{r}-1};1).\]
\end{theorem}

Viewed as a product of nested sums, the product $\zeta(\bk)\zeta(\bl)$ of two MZVs is expanded as a $\bbQ$-linear sum of MZVs. For example, when $k$ and $l$ are integers exceeding $1$, we obtain 
\[\zeta(k)\zeta(l) = \sum_{m,n>0}\frac{1}{m^{k}n^{l}} = \left(\sum_{0<m<n}+\sum_{0<n<m}+\sum_{0<m=n}\right)\frac{1}{m^{k}n^{l}}=\zeta(k,l)+\zeta(l,k)+\zeta(k+l).\]
This procedure, which is the so-called \emph{harmonic product} is naturally extended to arbitrary MZVs, as well as its iterated integral version (\emph{shuffle product}):
\begin{align}
    \zeta(2)^{2}
    &=\left(\int_{0<s_{1}<s_{2}<1}\frac{ds_{1}}{1-s_{1}}\frac{ds_{2}}{s_{2}}\right)\left(\int_{0<t_{1}<t_{2}<1}\frac{dt_{1}}{1-t_{1}}\frac{dt_{2}}{t_{2}}\right)\\
    &=\sum_{\substack{\bullet\in\{(s_{1}<s_{2}<t_{1}<t_{2})\\ (s_{1}<t_{1}<s_{2}<t_{2})\\ (t_{1}<s_{1}<s_{2}<t_{2})\\ (t_{1}<t_{2}<s_{1}<s_{2})\\ (t_{1}<s_{1}<t_{2}<s_{2})\\ (s_{1}<t_{1}<t_{2}<s_{2})\}}}\int_{0<\bullet<1}\frac{ds_{1}}{1-s_{1}}\frac{ds_{2}}{s_{2}}\frac{dt_{1}}{1-t_{1}}\frac{dt_{2}}{t_{2}}\\
    &=4\zeta(1,3)+2\zeta(2,2).
\end{align}

Comparing these two expansions, we obtain a big family of linear relations among MZVs. After some arguments extending it to the divergent cases, the family is the so-called \emph{regularized double shuffle relations}.
A modern tool handling it is non-commutative generating functions, as described in the next subsection.

\subsection{Non-commutative generating functions}
Hereafter, $R$ denotes a commutative $\bbQ$-algebra. Though there will appear several algebraic structures based on $R$, we omit the symbol $R$ unless needed. 

Let $\ff\coloneqq\ff_{\{e_{0},e_{1}\}}$ be the free Lie $R$-algebra generated by the formal symbols $e_{0}$ and $e_{1}$. By putting $\deg(e_{i})=1$, the Lie algebra $\ff$ becomes graded, and we can define the completion $\hff$ with respect to this grading. As is well-known, the completed universal enveloping algebra $\hcV\coloneqq U(\ff)^{\wedge}$, which is naturally isomorphic to the non-commutative power series algebra $R\ncjump{e_{0},e_{1}}$, has the natural structure as a completed Hopf $R$-algebra whose coproduct is defined by
\[\Delta(e_{i})\coloneqq 1\otimes e_{i}+e_{i}\otimes 1,\]
and its subset of all primitive elements is identified with $\hff$:
\[\hff = \{\varphi\in\hcV \mid \Delta(\varphi)=\varphi\otimes 1+1\otimes\varphi\}.\]
Similarly, the corresponding group $\exp(\hff)$ is identified with the subset of $\hcV$ consisting of all group-like elements:
\[\exp(\hff) = \{\varphi\in\hcV \mid \Delta(\varphi)=\varphi\otimes\varphi,~\varphi(0,0)=1\}.\]

For later convenience, we put
\[\fG \coloneqq \Ker(\hff\twoheadrightarrow\hff/[\hff,\hff]) = \{\varphi\in\hff\mid \varphi \equiv 0~(\mod~\deg\ge 2) \}\]
and
\[\cG \coloneqq \exp(\fG).\]

Moreover, denote by $\langle \mid \rangle$ the natural pairing on $\hcV\otimes\bbQ\langle e_{0},e_{1}\rangle$ extracting a coefficient: namely, for $\varphi\in\hcV$, we put
\[\varphi(e_{0},e_{1}) \eqqcolon \sum_{k=0}^{\infty}\sum_{a_{1},\ldots,a_{k}\in\{0,1\}}\langle \varphi\mid e_{a_{1}}\cdots e_{a_{k}}\rangle e_{a_{1}}\cdots e_{a_{k}},\]
and regard $\langle\varphi\mid\hyp\rangle$ as a $\bbQ$-linear map on the non-commutative polynomial algebra $\bbQ\langle e_{0},e_{1}\rangle$.

\begin{remark}\label{rem:shuffle_regularization}
Note that, as explained in \cite[\S 3.2]{furusho03}, an element $\varphi\in\cG$ is uniquely determined from the values $\{\langle\varphi \mid e_{0}w e_{1}\rangle\}_{w}$ given for each word $w\in\{e_{a_{1}}\cdots e_{a_{k}}\mid a_{i}\in\{0,1\}\}$.
\end{remark}

\subsection{Double shuffle set}
The subspace $\hcW\coloneqq \bbQ\oplus\hcV e_{1}$ has another (completed) Hopf algebra structure. Since $\hcW$ is isomorphic to the non-commutative power series algebra generated by $\{y_{k}\}_{k\ge 1} \coloneqq \{-e_{0}^{k-1}e_{1}\}_{k\ge 1}$ as a filtered\footnote{As the definition $y_{k}=-e_{0}^{k-1}e_{1}$ says, the filtration is defined by $\deg(y_{k})\coloneqq k$.} $R$-algebra, the map
\[\Delta^{\cW}(y_{k}) \coloneqq 1\otimes y_{k}+y_{k}\otimes 1+\sum_{i=1}^{k-1}y_{i}\otimes y_{k-i}\]
extends to the $R$-algebra morphism $\Delta^{\cW}$ on $\hcW$. This coproduct and counit $y_{k}\mapsto 0$ defines a completed Hopf algebra structure on $\hcW$. 

For each $\varphi\in\cG$, we put
\[\Gamma_{\varphi}(t) \coloneqq \exp\left(\sum_{k=2}^{\infty}\frac{(-1)^{k+1}}{k}\langle\varphi\mid e_{0}^{k-1}e_{1}\rangle t^{k}\right)\in R\jump{t},\]
and
\[\varphi_{\Ra}^{\ast}\coloneqq \Gamma_{\varphi}(-e_{1})^{-1}(1+\varphi_{1}e_{1}) \in \hcW,\]
where $\varphi_{1}e_{1}\in\hcW$ is the unique element determined by $\varphi = 1+\varphi_{0}e_{0}+\varphi_{1}e_{1}$.

\begin{definition}[Double shuffle set]
    For $\mu\in R$, we define
    \[\DMR_{\mu} \coloneqq \{\varphi\in\cG\mid \Delta^{\cW}(\varphi_{\Ra}^{\ast}) = \varphi_{\Ra}^{\ast}\otimes\varphi_{\Ra}^{\ast},~\langle\varphi\mid e_{0}e_{1}\rangle = \mu^{2}/24\}.\]
\end{definition}

Combining Remark \ref{rem:shuffle_regularization} and the shuffle product formula for iterated integrals, we see that the condition
\[\langle \Phi_{\KZ}\mid e_{0}^{k_{r}-1}e_{1}\cdots e_{0}^{k_{1}-1}e_{1}\rangle = (-1)^{r}\zeta(k_{1},\ldots,k_{r})\]
given for each positive integers $r,k_{1},\ldots,k_{r}$ ($k_{r}\ge 2$) determines the unique element $\Phi_{\KZ}\in\cG$. Then the regularized double shuffle relations for multiple zeta values is stated as follows.

\begin{theorem}[Regularized double shuffle relations for MZVs]
    We have $\Phi_{\KZ}\in\DMR_{2\pi i}$ for $R=\bbC$.
\end{theorem}

By passing to the dual, the defining properties of $\DMR_{\mu}$ are also described in terms of coefficients. This reformulation needs the precise definition of the shuffle and harmonic products. First, the shuffle product $\sh$ is defined to be the $\bbQ$-bilinear binary operation on $\bbQ\langle e_{0},e_{1}\rangle$, determined by the initial condition
\[w\sh 1=1\sh w=w\]
and the recursive rule
\[e_{i}w\sh e_{j}w' = e_{i}(w\sh e_{j}w')+e_{j}(e_{i}w\sh w'),\]
available for every $w,w_{1},w_{2}\in\bbQ\langle e_{0},e_{1}\rangle$. Next, the harmonic product $\ast$ is defined to be the $\bbQ$-bilinear binary operation on $\bbQ\langle y_{1},y_{2},\ldots\rangle$, determined by the initial condition
\[w\ast 1=1\ast w=w\]
and the recursive rule
\[y_{k}w\ast y_{l}w' = y_{k}(w\ast y_{l}w')+y_{l}(y_{k}w\ast w')+y_{k+l}(w\ast w'),\]
available for every $w,w_{1},w_{2}\in\bbQ\langle e_{0},e_{1}\rangle$. One may note that, the shuffle product always preserves both the degree and the depth ($\coloneqq$ $e_{1}$-degree), while the harmonic product only does for the degree. 

\begin{proposition}[{cf.~\cite{ikz06}, \cite{by24}}]\label{prop:rds}
    The set $\DMR_{\mu}$ consists of $\varphi\in\hcV$ satisfying:
    \begin{enumerate}
        \item $\varphi\equiv 1~(\mod~\deg\ge 2)$.
        \item $\langle\varphi\mid e_{0}e_{1}\rangle=\mu^{2}/24$.
        \item $\langle\varphi\mid w\sh w'\rangle=\langle\varphi\mid w\rangle\cdot\langle\varphi\mid w'\rangle$ for every $w,w'\in\bbQ\langle e_{0},e_{1}\rangle$.
        \item $\langle\varphi_{\Ra}^{\ast}\mid w\ast w'\rangle=\langle\varphi_{\Ra}^{\ast}\mid w\rangle\cdot\langle\varphi_{\Ra}^{\ast}\mid w'\rangle$ for every $w,w'\in\bbQ\langle y_{1},y_{2},\ldots\rangle$.
    \end{enumerate} 
\end{proposition}

\subsection{Recent embeddings and our main theorem}
The Kashiwara--Vergne conjecture also involves the solution set
\[
    \SolKV_{\mu}\subseteq\Aut(\hff).
\]
We use the injection
\begin{align}
    \aut_{-}^{\cV,(1)}\colon\cG &\to\Aut(\hff)\\
    \varphi &\mapsto\begin{cases}
        e_{0}\mapsto\varphi e_{0}\varphi^{-1},\\
        e_{1}\mapsto e_{1},
    \end{cases}
\end{align}
the image of $\DMR_{\mu}$ is conjectured to be $\SolKV_{\mu}$. Recently, Schneps \cite{schneps25}, Enriquez--Furusho \cite{ef25} and the first author \cite{kawamura25} independently showed that the image of $\DMR_{0}$ contained in $\KRV_{2}\coloneqq\SolKV_{0}$. In this paper, enhancing this result, we show the following theorem.

\begin{theorem}[$=$ Theorem \ref{thm:main_fourth_honbun}]\label{thm:main_fourth}
    Let $\varphi\in\cG$. Then the automorphism $\aut_{\varphi}^{\cV,(1)}$ solves the Kashiwara--Vergne conjecture if and only if, $\varphi$ satisfies the regularized double shuffle relations, in which the harmonic product formula is restricted for the case where one side is always depth $1$. Namely, defining
    \[\WDMR_{\mu}\coloneqq \{\varphi\colon\text{satisfies (1), (2), (3) in Prop.~\ref{prop:rds} and (4') below}\}\]
    with the instead condition
    \begin{enumerate}
        \item[(4')] $\langle\varphi_{\Ra}^{\ast}\mid y_{k}\ast w\rangle=\langle\varphi_{\Ra}^{\ast}\mid y_{k}\rangle\cdot\langle\varphi_{\Ra}^{\ast}\mid w\rangle$ for every $w\in\bbQ\langle y_{1},y_{2},\ldots\rangle$ and $k\ge 1$,
    \end{enumerate}
    we have $\WDMR_{\mu}=\SolKV_{\mu}$.
\end{theorem}

\begin{remark}\label{rem:ren_kv}
    Alekseev--Torossian \cite{at12} defined a symmetric subgroup $\KRV_{2}^{\mathrm{sym}}$ of the \emph{Kashiwara--Vergne group}. Kuno \cite{kuno25} interpreted it as an \emph{emergent} analogue of the Grothendieck--Teichm\"{u}ller group. Howarth and Ren \cite[Theorem C]{hr25} recently obtained a result for $\KRV_{2}^{\mathrm{sym}}$ similar to Theorem \ref{thm:main_fourth}, using Kuno's emergent pentagon equation. Their formulation uses the pentagon defect, which can be removed by \cite[Proposition 5.4]{hr25}.
\end{remark}

\subsection{Important lemmas, corollaries and applications}
Here we show the organization of this paper. In the first part, we investigate the so-called inertia-preserving condition in detail. In particular, we show that it is interpreted as a formalization of Bouillot's reduction of multitangent functions.
\begin{theorem}[$=$ Theorem \ref{thm:main_first_honbun}]\label{thm:main_first}
    Let us denote by $\cG^{\mu}_{\Red}$ the set consisting of elements satisfying Bouillot's reduction relation. Then we have $\cG^{\mu}_{\Red}=\cG^{\mu}_{\inert}$, where the right-hand side is the set satisfying inertia-preserving condition \cite[Definition 0.15]{ef25}.
\end{theorem}
The second part then studies interpretation of Kashiwara--Vergne structures as a part of double shuffle structures, with help of Ecalle's flexion theory. This part has Theorem \ref{thm:main_fourth} as its main content. In the final part, we consider about three involutions, coming from transpositions on $\{0,1,\infty\}$, naturally acting on $\cG^{\mu}_{\Red}=\cG^{\mu}_{\inert}$. Our result on this part is a multitangents-familiar interpretation of Enriquez--Furusho's involution $\Theta^{\mu}$, which also gives another proof of Ren's theorem \cite[Theorem 1.2]{ren26}. A new criterion for the so-called ``duality deduction'' is also shown:
\begin{theorem}[$=$ Theorem \ref{thm:main_extra_honbun}]\label{thm:main_extra}
    If the weak double shuffle relations (in the sense of Theorem \ref{thm:main_fourth}) imply the ``harmonic antipode relation'', the regularized double shuffle relations imply the duality relation.
\end{theorem}

\part{Multitangent functions and inertia-preserving condition}\label{part:multitangents}
\section{Some notions on generating functions and automorphisms}
\subsection{Notations}
Recall that
\begin{itemize}
    \item $R$ is a commutative $\bbQ$-algebra.
    \item $\ff\coloneqq\Lie_{R}(e_{0},e_{1})$ is the free Lie $R$-algebra generated by formal symbols $\{e_{0},e_{1}\}$,
    \item $\hff$ is the degree completion of $\ff$ by $\deg(e_{i})=1$,
    \item $\hcV\simeq R\ncjump{e_{0},e_{1}}$ is the degree completion of the universal enveloping algebra of $\ff$,
    \item $\fG$ is the degree $\ge 2$ part of $\hff$,
    \item $\cG\coloneqq\exp(\fG)$ is the corresponding group to $\fG$,
    \item $y_{k}\coloneqq -e_{0}^{k-1}e_{1}$ is an algebra generator of the below $\hcW$,
    \item $\Gamma_{\varphi}(t)\coloneqq\exp\left(\sum_{k\ge 1}(-1)^{k+1}\langle\varphi\mid e_{0}^{k-1}e_{1}\rangle t^{k}/k\right)$ is the Gamma function associated with $\varphi\in\hcV$, which is used as the correction term in the double shuffle structure,
\end{itemize}
and that $\hcW\coloneqq\bbQ\oplus\hcV e_{1}\simeq\bbQ\ncjump{y_{1},y_{2},\ldots}$ is a completed Hopf algebra with the coproduct
\[\Delta^{\cW}(y_{k})\coloneqq y_{k}\otimes 1+1\otimes y_{k}+\sum_{i=1}^{k-1}y_{i}\otimes y_{k-i}.\]
We also use the notation $\rev_{\re}$ (resp.~$\rev_{\ry}$) to indicate the reversing operator for words in $\{e_{0},e_{1}\}$ (resp.~$\{y_{k}\}_{k\ge 1}$). For example, we have
\[\rev_{\re}(e_{0}^{k-1}e_{1}e_{0}^{l-1}e_{1})=e_{1}e_{0}^{l-1}e_{1}e_{0}^{k-1},\]
and
\[\rev_{\ry}(y_{k}y_{l}) = y_{l}y_{k} = e_{0}^{l-1}e_{1}e_{0}^{k-1}e_{1}.\]
We note that, though $\rev_{\re}$ does not coincide with $\rev_{\ry}$ in general, one shows
\[e_{1}\rev_{\ry}(w) = \rev_{\re}(w)e_{1}\]
for all word in $\{y_{k}\}_{k\ge 1}$.
Moreover, it should be remarked that the antipode $S_{\sh}$ of the completed Hopf algebra $\hcV$ is given by
\[f\mapsto S_{\sh}(f)=\rev_{\re}(f(-e_{0},-e_{1})).\]

\subsection{Racinet product and Ihara bracket}\label{subsec:racinet}
We define the product $\varphi\circledast\psi$ by
\begin{align}
    (\varphi\circledast\psi)(e_{0},e_{1}) &\coloneqq \psi(\varphi(e_{0},e_{1})e_{0}\varphi(e_{0},e_{1})^{-1},e_{1})\varphi(e_{0},e_{1})\\
    &= \aut_{\varphi}^{\cV,(1)}(\psi)\varphi
\end{align}
for $\varphi,\psi\in\cG$. The corresponding bracket is defined as
\[\langle \varphi,\psi\rangle \coloneqq d_{\varphi}(\psi)-d_{\psi}(\varphi)+[\varphi,\psi]\]
for $\varphi,\psi\in\fG$. Here $d_{a}$ denotes the derivation
\[d_{a}\colon\fG\to\fG,\quad\begin{cases}
    e_{0} \mapsto 0,\\
    e_{1} \mapsto [e_{1},a].
\end{cases}\]

Here we note that $\aut_{\bullet}^{\cV,(1)}$ is \emph{not} compatible to the derivation $d_{\bullet}$, since our convention used to define $d_{\bullet}$ is chosen to be compatible with \cite{ef25}, not our $\aut_{\bullet}^{\cV,(1)}$.

\begin{lemma}[{\cite[Proposition 2.1 (i), (iv), (2.16), (2.14)]{racinet}}]\label{lem:MT}
    \begin{enumerate}
        \item The operation $\circledast$ is a group law on $\cG$.
        \item The correspondence
        \[
            R\longmapsto (\cG,\circledast)
        \]
        defines a pro-unipotent group. Its Lie algebra is $(\fG,\langle,\rangle)$.
        \item The Lie-exponential map $\fG\to\cG$ is given by
        \[\exp_{\circledast}(\varphi) \coloneqq \exp(s_{\varphi})(1),\]
        where $s_{\varphi}$ is the linear map sending $\psi\in\hcV$ to $\varphi\psi+d_{\varphi}(\psi)$.
    \end{enumerate}
\end{lemma}

Especially, the first assertion above gives a trivial torsor structure on $(\cG,\circledast)$.

\subsection{Harmonic regularization}
We also prepare some notations for coefficients.
For each $\varphi\in\hcV$, we consider the decomposition
\begin{align}
    \varphi &= \langle\varphi\mid 1\rangle + \varphi_{1}e_{1} + \varphi_{0}e_{0}\\ &= \langle\varphi\mid 1\rangle + e_{1}\varphi^{1} + e_{0}\varphi^{0}.
\end{align}
\begin{lemma}\label{lem:decomposition_formulas}
    For an element $\varphi\in\cG$, we have
    \begin{align}
        (\varphi^{-1})_{0} &= -\varphi^{-1}\varphi_{0}=-S_{\sh}(\varphi^{0}),\\
        (\varphi^{-1})_{1} &= -\varphi^{-1}\varphi_{1}=-S_{\sh}(\varphi^{1}),\\
        (\varphi^{-1})^{0} &= -\varphi^{0}\varphi^{-1}=-S_{\sh}(\varphi_{0}),\\
        (\varphi^{-1})^{1} &= -\varphi^{1}\varphi^{-1}=-S_{\sh}(\varphi_{1}).
    \end{align}
\end{lemma}
\begin{proof}
    For the first and second identity, we use $S_{\sh}(\varphi)=\varphi^{-1}$ and the formula 
    \begin{align}
    1 + (\varphi^{-1})_{0}e_{0} + (\varphi^{-1})_{1}e_{1}
    &= \varphi^{-1}\\
    &= \sum_{n=0}^{\infty}(-\varphi_{0}e_{0}-\varphi_{1}e_{1})^{n}\\
    &= 1+\left(\sum_{n=1}^{\infty} (-\varphi_{0}e_{0}-\varphi_{1}e_{1})^{n-1}\right)(-\varphi_{0}e_{0}-\varphi_{1}e_{1})\\
    &= 1+\varphi^{-1}(-\varphi_{0}e_{0}-\varphi_{1}e_{1}).
    \end{align}
    The rest two equalities are also shown in the same way.
\end{proof}
As seen in Introduction, the set $\DMR_{\mu}$ is defined to be the set of $\varphi\in\cG$ satisfies
\begin{equation}\label{eq:double_shuffle_equation}
    \Delta^{\cW}(\varphi^{\ast}_{\Ra}) = \varphi^{\ast}_{\Ra}\otimes\varphi^{\ast}_{\Ra},
\end{equation}
when putting $\varphi^{\ast}_{\Ra}\coloneqq\Gamma_{\varphi}(-e_{1})^{-1}(1+\varphi_{1}e_{1})$.
This condition is rephrased in terms of coefficients as follows: for $\varphi\in\hcV$, we define a $\bbQ$-linear map $\zeta^{\sh}_{\varphi}$ to be $\langle\varphi\mid -\rangle\circ\rev_{\re}$, namely
\[\varphi = \sum_{k=0}^{\infty}\sum_{a_{1},\ldots,a_{k}\in\{0,1\}}\zeta_{\varphi}^{\sh}(e_{a_{1}}\cdots e_{a_{n}})e_{a_{n}}\cdots e_{a_{1}}.\]
Moreover we put
\[\varphi^{\ast}_{\Ra} = \sum_{r=0}^{\infty}\sum_{k_{1},\ldots,k_{r}\ge 1}\zeta_{\varphi}^{\ast}(y_{k_{1}}\cdots y_{k_{r}})y_{k_{r}}\cdots y_{k_{1}}.\]
Then the condition $\varphi\in\cG$ exactly means that the map
\[\zeta^{\sh}_{\varphi}\colon\fH \coloneqq \bbQ\langle e_{0},e_{1}\rangle\to R\]
sends $1$ to $1$ and is a $\bbQ$-algebra homomorphism with respect to the shuffle product $\sh$. Moreover, under this, the condition \eqref{eq:double_shuffle_equation} amounts to say that
\[\zeta^{\ast}_{\varphi}\colon\fH^{1} \coloneqq \bbQ\langle y_{1},y_{2},\ldots\rangle \to R\]
is a $\bbQ$-algebra homomorphism with respect to the harmonic product $\ast$.

\subsection{Hurwitz-type values and coefficients}
A useful tool is an extension of the map $\zeta^{\ast}_{\varphi}$ on $\fH$.
Putting
\[\varphi_{\HK}^{\ast}\coloneqq\Gamma_{\varphi}(-e_{1})^{-1}\varphi\in\hcV,\]
we see that
\[1+(\varphi_{\HK}^{\ast})_{1}e_{1} = \varphi_{\Ra}^{\ast}\]
by definition, and
\[\varphi_{\HK}^{\ast} = \sum_{n=0}^{\infty}\left(\sum_{r=0}^{\infty}\sum_{k_{1},\ldots,k_{r}\ge 1}\zeta_{\varphi,n}^{\ast}(y_{k_{1}}\cdots y_{k_{r}})y_{k_{r}}\cdots y_{k_{1}}\right)e_{0}^{n},\]
where $\zeta_{\varphi,n}^{\ast}$ is the generalized value of $\zeta^{\ast}_{\varphi}$ defined by
\[\zeta_{\varphi,n}^{\ast}(y_{k_{1}}\cdots y_{k_{r}}) \coloneqq (-1)^{n}\sum_{\substack{n_{1},\ldots,n_{r}\ge 0\\ n_{1}+\cdots+n_{r}=n}}\left(\prod_{i=1}^{r}\binom{k_{i}+n_{i}-1}{n_{i}}\right)\zeta^{\ast}_{\varphi}(y_{k_{1}+n_{1}}\cdots y_{k_{r}+n_{r}}),\]
in the same argument as done in \cite[Proposition 2.9]{hk25}.
Related to it, this sequence $\{\zeta_{\varphi,n}^{\ast}(y_{k_{1}}\cdots y_{k_{r}})\}_{n}$ forms the MZV of \emph{Hurwitz type},
\[\zeta_{\varphi}^{\ast}(y_{k_{1}}\cdots y_{k_{r}};z)
    \coloneqq \sum_{n=0}^{\infty}\zeta_{\varphi,n}^{\ast}(y_{k_{1}}\cdots y_{k_{r}})z^{n},\]
as its generating function. Moreover, the generating function of these values with respect to $(k_{1},\ldots,k_{r})$ is expressed as follows:
\begin{lemma}\label{lem:hurwitz_gf}
    Let $\varphi\in\cG$. Then, the generating function
    \[\langle\varphi^{\ast}(e_{0},e_{1};z)\mid z^{n}\rangle = \sum_{r=0}^{\infty}\sum_{k_{1},\ldots,k_{r}\ge 1}\zeta_{\varphi,n}^{\ast}(y_{k_{1}}\cdots y_{k_{r}})y_{k_{r}}\cdots y_{k_{1}},\]
    which is associated with each $n$, has the expression
    \[\varphi_{\HK}^{\ast}(e_{0},e_{1}) = \sum_{n=0}^{\infty}\langle\varphi^{\ast}(e_{0},e_{1};z)\mid z^{n}\rangle e_{0}^{n}.\]
\end{lemma}
\begin{proof}
    For each $n\ge 0$, we define an element $A_{n}\in\hcW$ by
\[ \varphi_{\HK}^{\ast} = \sum_{n=0}^{\infty}A_{n}e_{0}^{n}.\]
Then the formula
\[A_n = \sum_{r\ge 1}\sum_{k_1,\ldots,k_r\ge 1}\zeta_{\varphi,n}^{\ast}(y_{k_1}\cdots y_{k_r})y_{k_{r}}\cdots y_{k_1}\]
is obvious. On the other hand, the definion of $\varphi^{\ast}(e_0,e_{1};z)$ yields
\begin{align}
\varphi^{\ast}(e_0,e_1;z)
&= \sum_{r\ge 0}\sum_{k_{1},\ldots,k_{r}\ge 1} \zeta_{\varphi}^{\ast}(y_{k_{1}}\cdots y_{k_{r}};z)y_{k_{r}}\cdots y_{k_{1}}\\
&= \sum_{n=0}^{\infty}A_nz^{n}.
\end{align}
Hence we obtain the conclusion.
\end{proof}

\begin{lemma}\label{lem:harmonic_hurwitz}
    When $\varphi\in\cG$, the generating function $\varphi^{\ast}(e_{0},e_{1};z)$ has the expression
    \[\varphi^{\ast}(e_{0},e_{1};z) = \varphi^{\ast}_{\Ra}(e_{0}-z,e_{1};0)\]
    as a $R\jump{z}\ncjump{e_{0},e_{1}}$.
\end{lemma}
\begin{proof}
    The right-hand side is computed as
    \begin{align}
        \varphi^{\ast}_{\Ra}(e_{0}-z,e_{1};0)
        &=\sum_{r=0}^{\infty}\sum_{k_{1},\ldots,k_{r}\ge 1}(-1)^{r}\zeta_{\varphi}^{\ast}(y_{k_{1}}\cdots y_{k_{r}})(e_{0}-z)^{k_{r}-1}e_{1}\cdots (e_{0}-z)^{k_{1}-1}e_{1}\\
        &=1+\sum_{r=1}^{\infty}\sum_{\substack{k_{1},\ldots,k_{r}\ge 1\\ 0\le c_{1}<k_{1}\\ \cdots\\ 0\le c_{r}<k_{r}}}(-1)^{r}
        \left(\prod_{j=1}^{r}\binom{k_{j}-1}{c_{j}}\right)\zeta_{\varphi}^{\ast}(y_{k_{1}}\cdots y_{k_{r}})
        (-z)^{c_{1}+\cdots+c_{r}}\\
        &\qquad\cdot e_{0}^{k_{r}-c_{r}-1}e_{1}\cdots
        e_{0}^{k_{1}-c_{1}-1}e_{1}\\
        &=1+\sum_{r=1}^{\infty}\sum_{l_{1},\ldots,l_{r}\ge 1}(-1)^{r}\zeta_{\varphi}^{\ast}(y_{l_{1}}\cdots y_{l_{r}};z)e_{0}^{l_{r}-1}e_{1}\cdots e_{0}^{l_{1}-1}e_{1}.
    \end{align}
    The resulting series is nothing but $\varphi^{\ast}(e_{0},e_{1};z)$.
\end{proof}

\subsection{Euler's theorem}
Finally, we introduce a condition related to Euler's result.
For each $\mu\in R$, we define
\[\cG^{\Eul{\mu}}(R)\coloneqq\left\{\varphi\in\cG~\left|~\langle \varphi\mid e_{0}^{2k-1}e_{1}\rangle = \frac{B_{2k}\mu^{2k}}{2(2k)!}~(k\ge 1)\right.\right\},\]
where $B_{2k}$ is the Bernoulli number. We usually omit the input $(R)$ from this notation.
Under the shuffle relation $\Phi_{\KZ}\in\cG$ of multiple zeta values, Euler's famous theorem is equivalent to:

\begin{theorem}[{Euler \cite{euler50}}]
    We have $\Phi_{\KZ}\in\cG^{\Eul{2\pi i}}(\bbC)$.
\end{theorem}

\subsection{Inertia-preserving structure}
\begin{definition}[Scaled BCH series]
    For $\mu\in R$, denote by $\bch_{\mu}(X,Y)$ the $\mu$-scaled version of the Baker--Campbell--Hausdorff (BCH) series. Precisely it is defined as follows: due to the usual BCH series
    \[\bch(X,Y) \coloneqq \log(\exp(X)\exp(Y))\]
    converges in the $(X,Y)$-adically completed algebra $\bbQ\ncjump{X,Y}$ and has no constant term, we can define
    \[\bch_{t}(X,Y) \coloneqq t^{-1}\bch(tX,tY) \in \bbQ[t]\ncjump{X,Y}.\]
    Hence we can define $\bch_{\mu}(X,Y)$ for each $\mu\in R$ as an image of $\bch_{t}(X,Y)$ under the map $\bbQ[t]\ncjump{X,Y} \to R\ncjump{X,Y}$ induced by $t\mapsto\mu$.
\end{definition}
\begin{definition}[{\cite[Definition 0.15]{ef25}}]
    \begin{enumerate}
        \item For $\mu\in R$, we define a subset $\cG_{\inert}^{\mu}$ of $\cG$ to be the set of all $\varphi\in\cG$ such that, there exists $\psi\in\exp(\hff)$ satisfying \[\bch_{\mu}(\varphi e_{0}\varphi^{-1},\psi e_{\infty}\psi^{-1})+e_{1}=0.\]
        \item We define a $R$-submodule $\fG_{\inert}$ of $\fG$ to be the set of all $\varphi\in\fG$ such that, there exists $\psi\in\fG$ satisfying \[[\varphi,e_{0}]+[\psi,e_{\infty}]=0.\]
    \end{enumerate}
\end{definition}
\begin{lemma}[{\cite[Lemma 15.17 (b), Lemma 0.7 (c)]{ef25}}]
    \begin{enumerate}
        \item The pair
        \[
            (\cG_{\inert}^{0},\cG_{\inert}^{\mu})
        \]
        is a subtorsor of $(\cG,\cG)$ coming from the product $\circledast$.
        \item We have a Lie subalgebra structure on $\fG_{\inert}$ of $(\fG,\langle,\rangle)$.
        \item The group $(\cG_{\inert}^{0},\circledast)$ is the corresponding group to the Lie algebra $(\fG_{\inert},\langle,\rangle)$.
    \end{enumerate}
\end{lemma}
\begin{lemma}\label{lem:ef16}
    We have
    \[\fG_{\inert} = \{\varphi\in\fG\mid \varphi_0-\varphi_1=\varphi^0-\varphi^1\}.\]
\end{lemma}
\begin{proof}
    We use the same convention as \cite[\S 16]{ef25}. Denote by $\bbL^{\wedge}$ the degree completion of $\bbL$, the free Lie algebra generated by $\{x,y\}$ (compatible with the completion $\ff\to\hff$) and prepare the decomposition
    \[f = f_{x}x+f_{y}y = xf^{x}+yf^{y}\]
    for each $f\in\bbL^{\wedge}$. Combining \cite[Lemma 16.3 (b)]{ef25} and \cite[Theorem 2.1]{schneps12}, we have
    \[\bbL_{\inert} = \{f \in \bbL^{\wedge} \mid f_{y}-f_{x} = -f^{y}+f^{x}\}.\]
    Hence we obtain the result due to the isomorphism $i\colon \bbL_{\inert,>1} \stackrel{\simeq}{\longrightarrow} \fG_{\inert}$ of \cite[Lemma 16.4]{ef25}.
\end{proof}

\subsection{Special automorphisms and derivations}
This section is only devoted to recall the definition of Kashiwara--Vergne structures. Instead one can read \cite{at12}, \cite{aet11} or \cite{fk}, etc.\\

As in the proof of Lemma \ref{lem:ef16}, we denote by $\bbL^{\wedge}$ the completed Lie algebra $\bbL$ generated by the symbols $\{x,y\}$ with the degree rule $\deg(x)=\deg(y)=1$. We denote $\bbL_{>1}^{\wedge}\coloneqq\overline{[\bbL^{\wedge},\bbL^{\wedge}]}$ (the overline denotes the closure). Its completed universal enveloping algebra is naturally regarded as the $R$-algebra of the non-commutative power series:
\[\hcV_{\bbL} \coloneqq R\ncjump{x,y}.\]
\begin{definition}[Special automorphism]\label{def:saut}
    Let $A$ be an automorphism of $\bbL^{\wedge}$. We say that $A$ is \emph{special} if it preserves $x+y$ and locally inner at each generator, namely it satisfies the following three condition:
    \begin{enumerate}
        \item $A(x)=fxf^{-1}$ for some $f\in\exp(\bbL^{\wedge})$.
        \item $A(y)=gyg^{-1}$ for some $g\in\exp(\bbL^{\wedge})$.
        \item $A(x+y)=x+y$.
    \end{enumerate}
    We denote by $\SAut_{2}\subset\Aut(\bbL^{\wedge})$ the subset consisting of all special automorphism of $\bbL^{\wedge}$. 
\end{definition}
\begin{definition}[Special derivation]\label{def:sder}
    Let $D$ be a derivation of $\bbL^{\wedge}$. We say that $D$ is \emph{special} if it kills $x+y$ and locally inner at each generator, namely it satisfies the following three condition:
    \begin{enumerate}
        \item $D(x)=[f,x]$ for some $f\in\bbL^{\wedge}$.
        \item $D(y)=[g,y]$ for some $g\in\bbL^{\wedge}$.
        \item $D(x+y)=0$.
    \end{enumerate}
    We denote by $\sder_{2}$ the $R$-module consisting of all special derivations of $\bbL^{\wedge}$. 
\end{definition}
\begin{remark}
    Every special automorphism (resp.~derivation) extendsuniquely to an algebra automorphism (resp.~derivation) on $\hcV_{\bbL}$.
\end{remark}
\begin{lemma}\label{lem:saut_and_sder}
    \begin{enumerate}
        \item $\SAut_{2}$ forms a pro-algebraic group.
        \item $\sder_{2}$ forms a completed Lie $R$-algebra associated with the pro-algebraic group $\SAut_{2}$.
    \end{enumerate}
\end{lemma}
\begin{proof}
    \begin{enumerate}
        We only put the proof of (1) (the second one is similar and easy). Let $A_{1}$ and $A_{2}$ be two special automorphisms and put $A_{i}(y)=g_{i}yg_{i}^{-1}$ for $i=1,2$. Then we have 
        \begin{align}
            (A_{1}\circ A_{2})(y)
            &= A_{1}(g_{2}yg_{2}^{-1})\\
            &= A_{1}(g_{2})g_{1}yg_{1}^{-1}A_{1}(g_{2})^{-1}.
        \end{align}
        Hence we see that $(A_{1}\circ A_{2})(y)$ is some conjugation of $y$. The same argument works to show
        \[(A_{1}\circ A_{2})(x) = A_{1}(f_{2})f_{1}xf_{1}^{-1}A_{1}(f_{2})^{-1}.\]           
        The rest condition $A_{1}\circ A_{2}\colon x+y\mapsto x+y$ is obivous.
    \end{enumerate}
\end{proof}
\begin{lemma}[{cf.~\cite[Remark 3.3]{at12}}]\label{lem:degree_1}
    If the elements $f$ and $g$ in Definition \ref{def:saut} exist, they are unique in $\exp(\bbL_{>1})$. In the situation of Definition \ref{def:sder}, they are unique in $\bbL_{>1}$.
\end{lemma}
\begin{proof}
    Let $u\in\{x,y\}$. We can check that an element in $\exp(\bbL^{\wedge})$ (resp.~$\bbL^{\wedge}$) commuting with $u$ must be written as $\exp(au)$ (resp.~$au$) by some $a\in R$. Such an element should be $1$ (resp.~$0$) under our current assumption of degree.
\end{proof}
Hereafter we denote by $\SAut_{2,>1}$ (resp.~$\sder_{2,>1}$) the subgroup (resp.~Lie subalgebra) consisting of elements sending each generator ($x$ and $y$) to what do not have the degree $2$ term in $\hcV_{\bbL}$.
\begin{remark}\label{rem:deg1}
    We easily check that the derivation $t$ defined by $x\mapsto [y,x]$ and $y\mapsto [x,y]$ is a special derivation, as well as $\Ad_{\exp(x+y)}$ is a special automorphism of $\bbL$. They are removed to construct to $\SAut_{2,>1}$ and $\sder_{2,>1}$: indeed we see that $\sder_{2}$ is decomposed to the direct sum of $Rt$ with $\sder_{2,>1}$.
\end{remark}
From Lemma \ref{lem:degree_1}, we can define the group injection $\rN^{-1}\colon\SAut_{2,>1}\to\exp(\bbL^{\wedge}_{>1})$; $A\mapsto g$ in the same notation, as well as the $R$-linear injection $\nu^{-1}\colon\sder_{2,>1}\to\bbL^{\wedge}_{>1}$; $D\mapsto g$. Thanks to the specialness, their inverses on the images are given as the restriction of the map
\[\rN\colon\exp(\bbL^{\wedge}_{>1})\to\Aut(\bbL^{\wedge}); g\mapsto\begin{cases}
    y \mapsto gyg^{-1}, & {}\\
    x+y \mapsto x+y, & {}
\end{cases}\]
and
\[\nu\colon\bbL^{\wedge}_{>1}\to\mathrm{Der}(\bbL^{\wedge}); g\mapsto\begin{cases}
    y \mapsto [g,y], & {}\\
    x+y \mapsto 0. & {}
\end{cases}\]
\begin{proposition}\label{prop:racinet_composition}
    Denote by $i$ the isomorphism $\bbL^{\wedge}\to\hff$ determined by $x\mapsto e_{0}$ and $y\mapsto -e_{1}$. We use the same symbol for the corresponding group isomorphism $\exp(\bbL^{\wedge})\to\exp(\hff)$. Then we have the following.
    \begin{enumerate}
        \item $\rN\circ i$ gives an isomorphism $\cG_{\inert}^{0}\simeq\SAut_{2,>1}$ of groups.
        \item $\nu\circ i$ gives an isomorphism $\fG_{\inert}\simeq\sder_{2,>1}$ of Lie algebras.
    \end{enumerate}
\end{proposition}
\begin{proof}
    The second statement follows from \cite[\S 16]{ef25}.
    Exponentiating gives (1).
    For product compatibility, Lemma \ref{lem:saut_and_sder} gives
    \[ 
        \rN\bigl(i^{-1}(i(g_{1})\circledast i(g_{2}))\bigr)=A_{1}\circ A_{2}.
    \]
\end{proof}
Define $\tr_2$ to be the quotient of $\hcV_{\bbL}$ by the commutator bracket, namely $\tr_2$ is the completion of the space spanned by all cyclic words in $\{x,y\}$. The natural projection is denoted by $|\cdot|$.\\

We now define the \emph{divergence} of each special derivation as follows. Recall that each $f\in\bbL^{\wedge}$ is decomposed as
\[f = f_{x}x+f_{y}y.\]
\begin{definition}[{\cite[Proposition 3.20, \S 5]{at12}}]
    \begin{enumerate}
        \item For a special derivation $D$ of $\bbL^{\wedge}$, write
        \[
            x\longmapsto[f,x],\qquad y\longmapsto[g,y].
        \]
        Then we define
        \[
            \div(D)\coloneqq |f_{x}x+g_{y}y|.
        \]
        \item Define the map $j\colon\SAut_{2,>1} \to \tr_{2}$ by
        \[j(\exp(u)) = \frac{\exp(u)-1}{u}\left(\div(u)\right),\]
        for each $u\in\sder_{2,>1}$. Here we used the fact that the definition of $\sder_{2}$ induces a canonical and well-defined action on $\tr_{2}$, not only on $\hcV_{\bbL}$.
    \end{enumerate}
\end{definition}
\begin{definition}\label{def:kv}
    \begin{enumerate}
        \item We define the \emph{Kashiwara--Vergne Lie algebra} $\krv_{2}$ as
        \[\krv_{2} \coloneqq \nu^{-1}(\{D\in\sder_{2}\mid \exists a\in R\jump{s},~\div(D) = |a(x+y)-a(x)-a(y)|\}).\]
        \item We define the \emph{Kashiwara--Vergne group} $\KRV_{2}$ as
        \[\KRV_{2} \coloneqq \rN^{-1}(\{A\in\SAut_{2}\mid \exists a\in R\jump{s},~j(A) = |a(x+y)-a(x)-a(y)|\}).\]
    \end{enumerate}
\end{definition}
\begin{remark}
    We here note that two slight (but delicate) differences between Definition \ref{def:kv} the ordinary convention. The \begin{enumerate}
        \item We defined the Kashiwara--Vergne group (resp.~Lie algebra) to be a subset of $\hcV_{\bbL}$ by taking $\rN^{-1}$ (resp.~$\nu^{-1}$), which is not done in the ordinary definition.
        \item The ``fundamental chord diagram'' $t$ (see Remark \ref{rem:deg1}) is usually included in the Kashiwara--Vergne structure, which we ignored by making the restriction ``$>1$''.
    \end{enumerate}
\end{remark}
Similarly, we also introduce $\SAut^{\mu}_{2}$ to be the set of an tomorphism $A$ of $\bbL^{\wedge}$ satisfying the following:
\begin{enumerate}
  \item $A(x) = f\exp(\mu x)f^{-1}$ for some $f\in\exp(\bbL^{\wedge})$.
  \item $A(y) = g\exp(\mu y)g^{-1}$ for some $g\in\exp(\bbL^{\wedge})$.
  \item $A(\bch_{\mu}(x,y)) = x+y$.
\end{enumerate}
Then a similar fact to Proposition \ref{prop:racinet_composition} holds as
\[ \rN^{\mu}\circ i\colon\cG^{\mu}_{\inert} \stackrel{\simeq}{\longrightarrow} \SAut_{2,>1}^{\mu},\]
where $\rN^{\mu}\colon\exp(\bbL^{\wedge}_{>1}) \to \Aut(\bbL^{\wedge})$ is given by
\[ g \mapsto \begin{cases}
  y \mapsto gyg^{-1}, & {}\\
  \bch_{\mu}(x,y) \mapsto x+y & {}.
\end{cases}\]
Note that it is also shown in \cite[Lemma 14.11]{ef25}. Using this set, we can arrive at the definition of Kashiwara--Vergne torsor:
\begin{definition}
  We define the solution set of Kashiwara--Vergne problem to be
  \[\SolKV_{\mu} \coloneqq (\rN^{\mu})^{-1}(\{A\in\SAut_{2}^{\mu}\mid \exists a\in R\jump{s},~j(A) = |a(x+y)-a(x)-a(y)|\}).\]
\end{definition}

\section{Bouillot reduction relation}
In this section, we work with the fixed $\mu\in R$.
\begin{definition}
    For each $\varphi\in\cG$ and $k_{1},\ldots,k_{r}\ge 1$, we define
    \begin{align}
        \hzeta_{\varphi}(y_{k_{1}}\cdots y_{k_{r}};z)
        &\coloneqq\sum_{i=0}^{r}(-1)^{k_{1}+\cdots+k_{i}}
        \zeta_{\varphi}^{\ast}(y_{k_{i}}\cdots y_{k_{1}};-z)
        \zeta_{\varphi}^{\ast}(y_{k_{i+1}}\cdots y_{k_{r}};z)\nonumber\\
        &\quad+\sum_{i=1}^{r}(-1)^{k_{1}+\cdots+k_{i-1}}
        \zeta_{\varphi}^{\ast}(y_{k_{i-1}}\cdots y_{k_{1}};-z)
        \frac{1}{z^{k_{i}}}
        \zeta_{\varphi}^{\ast}(y_{k_{i+1}}\cdots y_{k_{r}};z).
    \end{align}
    as a formal Laurent series in $R(\!(z)\!)$. We always understand that $\hzeta_{\varphi}(1;z)=1$.
\end{definition}

This definition reads as follows in terms of generating functions.

\begin{lemma}\label{lem:expression_T}
    Let $\varphi\in\cG$. Then we have
    \begin{align}
        \cT_{\varphi}(e_{0},e_{1};z)
        &\coloneqq\sum_{r=0}^{\infty}\sum_{k_{1},\ldots,k_{r}\ge 1}\hzeta_{\varphi}(y_{k_{1}}\cdots y_{k_{r}};z)y_{k_{r}}\cdots y_{k_{1}}\\
        &=\varphi^{\ast}_{\Ra}(e_{0}-z,e_{1})\cdot\left(1-\frac{1}{z-e_{0}}e_{1}\right)\cdot\rev_{\ry}(\varphi^{\ast}_{\Ra}(z-e_{0},-e_{1})),
    \end{align}
    where the element
    \[\frac{1}{z-e_{0}}=\frac{1}{z}+\frac{1}{z^{2}}e_{0}+\frac{1}{z^{3}}e_{0}^{2}+\cdots\]
    makes sense as an element of $R(\!(z)\!)\ncjump{e_{0},e_{1}}$.
\end{lemma}
\begin{proof}
    First we see that
    \begin{align}
        \cT_{\varphi}(e_{0},e_{1};z)
        &=\sum_{0\le i\le r}\sum_{k_{1},\ldots,k_{r}\ge 1}
        (-1)^{k_{1}+\cdots+k_{i}}
        \zeta_{\varphi}^{\ast}(y_{k_{i}}\cdots y_{k_{1}};-z)
        \zeta_{\varphi}^{\ast}(y_{k_{i+1}}\cdots y_{k_{r}};z)
        y_{k_{r}}\cdots y_{k_{1}}\nonumber\\
        &\quad+\sum_{1\le i\le r}\sum_{k_{1},\ldots,k_{r}\ge 1}
        (-1)^{k_{1}+\cdots+k_{i-1}}
        \zeta_{\varphi}^{\ast}(y_{k_{i-1}}\cdots y_{k_{1}};-z)
        \frac{1}{z^{k_{i}}}\nonumber\\
        &\qquad\cdot\zeta_{\varphi}^{\ast}(y_{k_{i+1}}\cdots y_{k_{r}};z)
        y_{k_{r}}\cdots y_{k_{1}}\\
        &=\left(\sum_{r=0}^{\infty}\sum_{k_{1},\ldots,k_{r}\ge 1}
        \zeta_{\varphi}^{\ast}(y_{k_{1}}\cdots y_{k_{r}};z)y_{k_{r}}\cdots y_{k_{1}}\right)\nonumber\\
        &\qquad\cdot\left(\sum_{s=0}^{\infty}\sum_{l_{1},\ldots,l_{s}\ge 1}
        (-1)^{l_{1}+\cdots+l_{s}}\zeta_{\varphi}^{\ast}(y_{l_{1}}\cdots y_{l_{s}};-z)
        y_{l_{1}}\cdots y_{l_{s}}\right)\nonumber\\
        &\quad+\left(\sum_{r=0}^{\infty}\sum_{k_{1},\ldots,k_{r}\ge 1}
        \zeta_{\varphi}^{\ast}(y_{k_{1}}\cdots y_{k_{r}};z)y_{k_{r}}\cdots y_{k_{1}}\right)
        \left(\sum_{k\ge 1}z^{-k}y_{k}\right)\nonumber\\
        &\qquad\cdot\left(\sum_{s=0}^{\infty}\sum_{l_{1},\ldots,l_{s}\ge 1}
        (-1)^{l_{1}+\cdots+l_{s}}\zeta_{\varphi}^{\ast}(y_{l_{1}}\cdots y_{l_{s}};-z)
        y_{l_{1}}\cdots y_{l_{s}}\right)\\
        &=\varphi^{\ast}(e_{0},e_{1};z)
        \left(1-\frac{1}{z-e_{0}}e_{1}\right)
        \rev_{\ry}(\varphi^{\ast}(-e_{0},-e_{1};-z)).
    \end{align}
    Combining it with Lemma \ref{lem:harmonic_hurwitz}, we obtain the desired result.
\end{proof}

Next, we consider the formalization of the ``reduced'' multitangents.
First, we define
\[ \delta^{\mu}(y_{k_{1}}\cdots y_{k_{r}}) \coloneqq \begin{cases}
    (\mu/2)^{r}/r! & \text{if }r\colon\text{even and }\bk=\underbrace{1,\ldots,1}_{r},\\
    0 & \text{otherwise},
\end{cases}\]
\begin{align*}
    \cR_{\varphi}(e_{0},e_{1};z)
    &\coloneqq \sum_{r=0}^{\infty}\sum_{k_{1},\ldots,k_{r}}
    \delta^{\mu}(y_{k_{1}}\cdots y_{k_{r}})y_{k_{r}}\cdots y_{k_{1}}\\
    &\quad+\sum_{r=0}^{\infty}\sum_{k_{1},\ldots,k_{r}}
    \sum_{i=1}^{r}\sum_{\substack{a+j+b=k_{i}\\ a,b\ge 0\\ j\ge 1}}
    (-1)^{k_{1}+\cdots+k_{i-1}+a}
    \zeta_{\varphi,a}^{\ast}(y_{k_{i-1}}\cdots y_{k_{1}})\\
    &\qquad\cdot\hzeta_{\varphi}(y_j;z)
    \zeta_{\varphi,b}^{\ast}(y_{k_{i+1}}\cdots y_{k_{r}})
    y_{k_{r}}\cdots y_{k_{1}}.
\end{align*}
as an element of $R(\!(z)\!)\ncjump{e_{0},e_{1}}$.

\begin{remark}
    Here we briefly desribe the original motivation. Recalling the KZ-associator $\Phi_{\KZ}$ has the expression
    \[(\Phi_{\KZ})^{\ast}_{\Ra} = \sum_{r=0}^{\infty}\sum_{k_{1},\ldots,k_{r}\ge 1}\zeta^{\ast}(k_{1},\ldots,k_{r})y_{k_{r}}\cdots y_{k_{1}}\]
    via the harmonic-regularized multiple zeta values $\{\zeta^{\ast}(\bk)\}_{\bk}$, we rephrase Bouillot's reduction theorem in terms of generating functions as follows.
\end{remark}

\begin{theorem}[{Bouillot \cite[Theorem 6]{bouillot14}}]
    We have $\cR_{\Phi_{\KZ}}=\cT_{\Phi_{\KZ}}$.
\end{theorem}

To consider the same condition as Bouillot's original theorem, we should find a nice expression of $\cR_{\varphi}$ for a general $\varphi$.
The following assertion shall be a candidate; to state it, we introduce a modified generating function of Bernoulli numbers as
\[B_{\mu}(t) \coloneqq \frac{\mu}{2}\frac{e^{\mu t}+1}{e^{\mu t}-1}-\frac{1}{t} = \sum_{k=1}^{\infty}\frac{B_{2k}\mu^{2k}}{(2k)!}t^{2k-1}.\]
\begin{lemma}\label{lem:expression_R}
    Assume $\varphi\in\cG^{\Eul{\mu}}$. Then we have
    \begin{align*}
        \cR_{\varphi}(e_{0},e_{1};z)
        &=\cos\left(\frac{\mu i e_{1}}{2}\right)\\
        &\quad-\Gamma_{\varphi}(-e_{1})^{-1}\varphi
        \left(\frac{1}{z-e_{0}}-B_{\mu}(e_{0}-z)\right)
        \varphi^{-1}e_{1}\Gamma_{\varphi}(e_{1})^{-1}.
    \end{align*}
\end{lemma}
\begin{proof}
    Values of the constant part $\delta^{\mu}(-)$ are summed up as
    \begin{align}
        \sum_{r=0}^{\infty}\sum_{k_{1},\ldots,k_{r}\ge 1}\delta^{\mu}(y_{k_{1}}\cdots y_{k_{r}})y_{k_{r}}\cdots y_{k_{1}}
        &=\sum_{s=0}^{\infty}\frac{(\mu/2)^{2s}}{(2s)!}y_{1}^{2s}\\
        &=\cos(\mu ie_{1}/2).
    \end{align}
    For the rest part, we can do the similar argument to the proof of Lemma \ref{lem:expression_T} as
    \begin{align}
        &\sum_{r=0}^{\infty}\sum_{k_{1},\ldots,k_{r}\ge 1}
        \sum_{i=1}^{r}\sum_{\substack{a+s+b=k_{i}\\ a,b\ge 0\\ s\ge 1}}
        (-1)^{k_{1}+\cdots+k_{i-1}+a}
        \zeta_{\varphi,a}^{\ast}(y_{k_{i-1}}\cdots y_{k_{1}})\nonumber\\
        &\qquad\cdot\hzeta_{\varphi}(y_s;z)
        \zeta_{\varphi,b}^{\ast}(y_{k_{i+1}}\cdots y_{k_{r}})
        y_{k_{r}}\cdots y_{k_{1}}\\
        &=\sum_{\substack{a,b\ge 0\\ j\ge 1}}
        \left(\sum_{r\ge 0}\sum_{k_{1},\ldots,k_{r}\ge 1}
        \zeta^{\ast}_{\varphi,b}(y_{k_{1}}\cdots y_{k_{r}})
        y_{k_{r}}\cdots y_{k_{1}}e_{0}^{b}\right)\nonumber\\
        &\qquad\cdot\hzeta_{\varphi}(y_j;z)e_{0}^{j-1}
        \left(\sum_{s\ge 0}\sum_{l_{1},\ldots,l_{s}\ge 1}
        (-1)^{l_{1}+\cdots+l_{s}+a}
        \zeta^{\ast}_{\varphi,a}(y_{l_{1}}\cdots y_{l_{s}})
        e_{0}^{a}(-e_{1})y_{l_{1}}\cdots y_{l_{s}}\right)\\
        &=\varphi_{\HK}^{\ast}(e_{0},e_{1})
        \left(\sum_{j\ge 1}\hzeta_{\varphi}(y_j;z)e_{0}^{j-1}\right)
        \rev_{\ry}(\varphi_{\HK}^{\ast}(-e_{0},-e_{1})(-e_{1}))\\
        &=\Gamma_{\varphi}(-e_{1})^{-1}\varphi(e_{0},e_{1})
        \left(\sum_{j\ge 1}\hzeta_{\varphi}(y_j;z)e_{0}^{j-1}\right)\rev_{\re}(\varphi_{\HK}^{\ast}(-e_{0},-e_{1}))(-e_{1})\\
        &=-\Gamma_{\varphi}(-e_{1})^{-1}\varphi
        \left(\sum_{j\ge 1}\hzeta_{\varphi}(y_j;z)e_{0}^{j-1}\right)
        \varphi^{-1}e_{1}\Gamma_{\varphi}(e_{1})^{-1},
    \end{align}
    where we used Lemma \ref{lem:hurwitz_gf} in the second equality. Hence it suffices to show
    \[\sum_{j\ge 1}\hzeta_{\varphi}(y_j;z)e_{0}^{j-1} = \frac{1}{z-e_{0}} - B_{\mu}(e_{0}-z).\]
    Due to the assumption $\varphi\in\cG^{\Eul{\mu}}$, we have
    \begin{align}
        &\sum_{j\ge 1}\hzeta_{\varphi}(y_j;z)e_{0}^{j-1}\\
        &=\sum_{j\ge 1}((-1)^{j}\zeta_{\varphi}^{\ast}(y_j;-z)+\zeta_{\varphi}^{\ast}(y_j;z)+z^{-j})e_{0}^{j-1}\\
        &=\frac{1}{z-e_{0}}+\sum_{j=1}^{\infty}\sum_{m=0}^{\infty}\binom{j+m-1}{m}\zeta_{\varphi}^{\ast}(y_{j+m})\left(-z^{m}(-e_{0})^{j-1}+(-z)^{m}e_{0}^{j-1}\right)\\
        &=\frac{1}{z-e_{0}}+\sum_{N=1}^{\infty}\zeta_{\varphi}^{\ast}(y_N)(-(z-e_{0})^{N-1}+(e_{0}-z)^{N-1})\\
        &=\frac{1}{z-e_{0}}+\sum_{N=1}^{\infty}2\zeta_{\varphi}^{\ast}(y_{2N})(e_{0}-z)^{2N-1}\\
        &=\frac{1}{z-e_{0}}-\sum_{N=1}^{\infty}\frac{B_{2N}}{(2N)!}\mu^{2N}(e_{0}-z)^{2N-1}\\
        &=\frac{1}{z-e_{0}} - B_{\mu}(e_{0}-z),
    \end{align}
    as desired.
\end{proof}

For our purpose, we consider the series
\[
    \cB_{\varphi}(e_{0},e_{1};z)
    \coloneqq \partial_{e_{1}}\left(
    \Gamma_{\varphi}(-e_{1})
    \left(\cR_{\varphi}(e_{0},e_{1};z)-\cT_{\varphi}(e_{0},e_{1};z)\right)
    \Gamma_{\varphi}(e_{1})\right)\varphi.
\]
Since both $\cT_{\varphi}$ and $\cR_{\varphi}$ belong to $R(\!(z)\!)\ncjump{y_{1},y_{2},\ldots}$ and have the same constant term, the factor
\[\Gamma_{\varphi}(-e_{1})\left(\cT_{\varphi}(e_{0},e_{1};z)-\cR_{\varphi}(e_{0},e_{1};z)\right)\Gamma_{\varphi}(e_{1})\]
is in $R(\!(z)\!)\ncjump{e_{0},e_{1}}e_{1}$ and thus the above $\partial_{e_{1}}$ eliminates no information. Hence Bouillot's result is equivalent to $\cB_{\Phi_{\KZ}}=0$.
The following is convenient properties for $\cB_{\varphi}$.

\begin{proposition}\label{prop:reduction_expression}
    Let $\varphi\in\cG^{\Eul{\mu}}$. Then we have the following.
    \begin{enumerate}
        \item $\cB_{\varphi}(e_{0},e_{1};z)$ has no polar term; namely it belongs to $R\jump{z}\ncjump{e_{0},e_{1}}$. 
        \item $\cB_{\varphi}(e_{0}-z,e_{1};0) = \cB_{\varphi}(e_{0},e_{1};z)$.
        \item When $z=0$, it has the symmetric expression
        \begin{align*}
            \cB_{\varphi}(e_{0},e_{1};0)
            &=\varphi(B_{\mu}(e_{0})-\varphi^{0})
            +(B_{\mu}(e_{1})-\varphi_{1})\varphi\\
            &\quad+\varphi_{0}(1+e_{0}\varphi^{0})
            +(1+\varphi_{1}e_{1})\varphi^{1}.
        \end{align*}
    \end{enumerate}
\end{proposition}
\begin{proof}
    By Lemma \ref{lem:expression_T}, we have $\cT_{\varphi}(e_{0}-z,e_{1};0)=\cT_{\varphi}(e_{0},e_{1};z)$. Hence we have to show
    \[\cR_{\varphi}(e_{0}-z,e_{1};0)=\cR_{\varphi}(e_{0},e_{1};z),\]
    which follows from the translation invariance $\varphi(e_{0}-z,e_{1})$ given by $\varphi\in\cG$. It shows (2).\\
    Next we show (1). Lemma \ref{lem:expression_T} yields 
    \begin{align}
        &\Gamma_{\varphi}(-e_{1})\cT_{\varphi}(e_{0},e_{1};z)\Gamma_{\varphi}(e_{1})\\
        &=(1+\varphi_{1}(e_{0}-z,e_{1})e_{1})\cdot\left(1-\frac{1}{z-e_{0}}e_{1}\right)\cdot\rev_{\ry}(1-\varphi_{1}(z-e_{0},-e_{1})e_{1})\\
        &=(1+\varphi_{1}(e_{0}-z,e_{1})e_{1})\cdot\left(1-\frac{1}{z-e_{0}}e_{1}\right)\cdot(1-\rev_{\re}(\varphi_{1}(z-e_{0},-e_{1})e_{1}))\\
        &=(1+\varphi_{1}(e_{0}-z,e_{1})e_{1})\cdot\left(1-\frac{1}{z-e_{0}}e_{1}\right)\cdot(1-S_{\sh}(\varphi_{1})(e_{0}-z,e_{1})e_{1})\\
        &=(1+\varphi_{1}(e_{0}-z,e_{1})e_{1})
        (1-S_{\sh}(\varphi_{1})(e_{0}-z,e_{1})e_{1})\nonumber\\
        &\quad-(1+\varphi_{1}(e_{0}-z,e_{1})e_{1})
        \frac{1}{z-e_{0}}
        (1-e_{1}S_{\sh}(\varphi_{1})(e_{0}-z,e_{1}))e_{1}.\\
    \end{align}
    On the other hand, Lemma \ref{lem:expression_R} shows
    \begin{align}
        &\Gamma_{\varphi}(-e_{1})\left(\cR_{\varphi}(e_{0},e_{1};z)-\cos\left(\frac{\mu i e_{1}}{2}\right)\right)\Gamma_{\varphi}(e_{1})\\
        &=-\varphi\left(\frac{1}{z-e_{0}}-B_{\mu}(e_{0}-z)\right)\varphi^{-1}e_{1}\\
        &=\varphi B_{\mu}(e_{0}-z)\varphi^{-1}e_{1}-\varphi(e_{0}-z,e_{1})\frac{1}{z-e_{0}}S_{\sh}(\varphi(e_{0}-z,e_{1}))e_{1}\\
        &=\varphi B_{\mu}(e_{0}-z)\varphi^{-1}e_{1}
        -(1+\varphi_{1}(e_{0}-z,e_{1})e_{1}
        +\varphi_{0}(e_{0}-z,e_{1})(e_{0}-z))\nonumber\\
        &\qquad\cdot\frac{1}{z-e_{0}}
        (1-(e_{0}-z)S_{\sh}(\varphi_{0})(e_{0}-z,e_{1})
        -e_{1}S_{\sh}(\varphi_{1})(e_{0}-z,e_{1}))e_{1}\\
        &=\varphi B_{\mu}(e_{0}-z)\varphi^{-1}e_{1}
        +\varphi_{0}(e_{0}-z,e_{1})\varphi^{-1}e_{1}
        -\varphi S_{\sh}(\varphi_{0})(e_{0}-z,e_{1})e_{1}\nonumber\\
        &\quad+\varphi_{0}(e_{0}-z,e_{1})(e_{0}-z)
        S_{\sh}(\varphi_{0})(e_{0}-z,e_{1})e_{1}\nonumber\\
        &\quad+\Gamma_{\varphi}(-e_{1})\cT_{\varphi}(e_{0},e_{1};z)
        \Gamma_{\varphi}(e_{1})\nonumber\\
        &\quad-(1+\varphi_{1}(e_{0}-z,e_{1})e_{1})
        (1-S_{\sh}(\varphi_{1})(e_{0}-z,e_{1})e_{1}).
    \end{align}
    Hence we obtain
    \begin{align}
        &\cB_{\varphi}(e_{0},e_{1};z)\\
        &=\partial_{e_{1}}\left(\Gamma_{\varphi}(-e_{1})\left(\cR_{\varphi}(e_{0},e_{1};z)-\cT_{\varphi}(e_{0},e_{1};z)\right)\Gamma_{\varphi}(e_{1})\right)\varphi\\
        &=\Gamma_{\varphi}(-e_{1})\Gamma_{\varphi}(e_{1})
        \cos\left(\frac{\mu i e_{1}}{2}\right)
        +\varphi B_{\mu}(e_{0}-z)+\varphi_{0}(e_{0}-z,e_{1})\nonumber\\
        &\quad-\varphi S_{\sh}(\varphi_{0})(e_{0}-z,e_{1})\varphi
        +\varphi_{0}(e_{0}-z,e_{1})(e_{0}-z)
        S_{\sh}(\varphi_{0})(e_{0}-z,e_{1})\varphi\nonumber\\
        &\quad-\partial_{e_{1}}\left(
        (1+\varphi_{1}(e_{0}-z,e_{1})e_{1})
        (1-S_{\sh}(\varphi_{1})(e_{0}-z,e_{1})e_{1})\right)\varphi,
    \end{align}
    which shows (1), since it is an well-defined element as the value of $\cB_{\varphi}$ at $z=0$. For (3), we continude the computation as
    \begin{align}
        &\cB_{\varphi}(e_{0},e_{1};0)\\
        &=\partial_{e_{1}}\left(\Gamma_{\varphi}(-e_{1})\Gamma_{\varphi}(e_{1})
        \cos\left(\frac{\mu i e_{1}}{2}\right)\right)\varphi
        +\varphi B_{\mu}(e_{0})+\varphi_{0}\nonumber\\
        &\quad-\varphi S_{\sh}(\varphi_{0})\varphi
        +\varphi_{0}e_{0}S_{\sh}(\varphi_{0})\varphi
        -\partial_{e_{1}}\left((1+\varphi_{1}e_{1})
        (1-S_{\sh}(\varphi_{1})e_{1})\right)\varphi\\
        &=\partial_{e_{1}}\left(\Gamma_{\varphi}(-e_{1})\Gamma_{\varphi}(e_{1})
        \cos\left(\frac{\mu i e_{1}}{2}\right)\right)\varphi
        +\varphi B_{\mu}(e_{0})+\varphi_{0}-\varphi\varphi^{0}\nonumber\\
        &\quad+\varphi_{0}e_{0}\varphi^{0}
        -\partial_{e_{1}}\left((1+\varphi_{1}e_{1})
        (1-\varphi^{1}\varphi^{-1}e_{1})\right)\varphi\\
        &=\partial_{e_{1}}\left(\Gamma_{\varphi}(-e_{1})\Gamma_{\varphi}(e_{1})
        \cos\left(\frac{\mu i e_{1}}{2}\right)\right)\varphi
        +\varphi B_{\mu}(e_{0})+\varphi_{0}-\varphi\varphi^{0}\nonumber\\
        &\quad+\varphi_{0}e_{0}\varphi^{0}-\varphi_{1}\varphi
        +\varphi^{1}+\varphi_{1}e_{1}\varphi^{1}.
    \end{align}    
    Then we can apply a standard argument
    \begin{align}
        \Gamma_{\varphi}(-e_{1})\Gamma_{\varphi}(e_{1})
        &=\exp\left(-\sum_{k=1}^{\infty}\frac{\zeta_{\varphi}(y_{2k})}{k}e_{1}^{2k}\right)\\
        &=\exp\left(-\sum_{k=1}^{\infty}\frac{B_{2k}(\mu e_{1})^{2k}}{2k(2k)!}\right)\\
        &=\frac{\mu i e_{1}/2}{\sin(\mu i e_{1}/2)}
    \end{align}
    and then have
    \[\Gamma_{\varphi}(-e_{1})\Gamma_{\varphi}(e_{1})\cos\left(\frac{\mu i e_{1}}{2}\right)=1+B_{\mu}(e_{1})e_{1},\]
    which completes the proof.
\end{proof}

Now we can define one of our main object $\cG^{\mu}_{\Red}$. It is the set of all element $\varphi$ satisfying the identity $\cR_{\varphi}=\cT_{\varphi}$ of Bouillot-type. It is reformulated as follows.

\begin{definition}
    We put
    \begin{align}
        \cG^{\mu}_{\Red}
        &=\{\varphi\in\cG^{\Eul{\mu}}\mid \cR_{\varphi}(e_{0},e_{1};z)=\cT_{\varphi}(e_{0},e_{1};z)\}\\
        &=\{\varphi\in\cG^{\Eul{\mu}}\mid \cB_{\varphi}(e_{0},e_{1};z)=0\}.
    \end{align}
    Thanks to Proposition \ref{prop:reduction_expression}, it is also written as
    \begin{equation}\label{eq:reduction_condition}
        \cG^{\mu}_{\Red}=\{\varphi\in\cG\mid \varphi(B_{\mu}(e_{0})-\varphi^{0})+(B_{\mu}(e_{1})-\varphi_{1})\varphi+\varphi_{0}(1+e_{0}\varphi^{0})+(1+\varphi_{1}e_{1})\varphi^{1}=0\}.
    \end{equation}
\end{definition}
\begin{remark}
    In the definition \eqref{eq:reduction_condition}, we do not need to require $\varphi\in\cG^{\Eul{\mu}}$. Indeed, due to $\varphi\in\cG$, we have
    \[\langle\varphi\mid e_{0}^{n}\rangle=\langle\varphi_{0}\mid e_{0}^{n}\rangle=\langle\varphi^{0}\mid e_{0}^{n}\rangle=0,\]
    for every $n\ge 1$. Hence the coefficient of the condition in \eqref{eq:reduction_condition} at $e_{0}^{2k-1}$ ($k\ge 1$) is
    \begin{align}
        0&= \langle\varphi(B_{\mu}(e_{0})-\varphi^{0})+(B_{\mu}(e_{1})-\varphi_{1})\varphi+\varphi_{0}(1+e_{0}\varphi^{0})+(1+\varphi_{1}e_{1})\varphi^{1}\mid e_{0}^{2k-1}\rangle\\
        &=\langle B_{\mu}(e_{0})-\varphi_{1}+\varphi^{1}\mid e_{0}^{2k-1}\rangle\\
        &=\frac{B_{2k}}{(2k)!}\mu^{2k}-\langle\varphi\mid e_{0}^{2k-1}e_{1}-e_{1}e_{0}^{2k-1}\rangle\\
        &=\frac{B_{2k}}{(2k)!}\mu^{2k}-2\langle\varphi\mid e_{0}^{2k-1}e_{1}\rangle,
    \end{align}
    which automatically deduces $\varphi\in\cG^{\Eul{\mu}}$.
\end{remark}

\section{Inertia-preserving condition and reduction relations}
Hereafter, we use the notation \[e_{\infty} \coloneqq -e_{0}-e_{1}.\]

\subsection{Proof of ``$\Red=\inert$''}
\begin{lemma}\label{lem:reduction_lie_algebra}
    The subspace
    \[\fG_{\Red} \coloneqq \{\varphi\in\fG\mid \varphi_0-\varphi_1=\varphi^0-\varphi^1\}\subseteq\fG\]
    is the Lie algebra associated with the pro-algebraic group $(\cG^{0}_{\Red},\circledast)$.
\end{lemma}
\begin{proof}
    It is due to a usual procedure of linearization.
\end{proof}

\begin{lemma}\label{lem:schneps_fact}
Assume $h\in\hff$. Then the following are equivalent:
\begin{enumerate}
\item There exist elements $P,Q\in \hcV$ such that $h=Pe_{\infty}Q$ and $P(0,0)=Q(0,0)=1$.
\item There exists an element $\psi\in\exp(\hff)$ such that $h=\psi e_{\infty}\psi^{-1}$.
\end{enumerate} 
\end{lemma}
\begin{proof}
The implication (2) $\implies$ (1) is obvious. We assume (1) to show the converse. Denote by $f\coloneqq h-e_{\infty}$; it does not have the constant term from the assumption $P(0,0)=Q(0,0)=1$. We define the algebra automorphism $c_{z}$ on $\hcV\jump{z}$ by $e_{0}\mapsto e_{0}$, $z\mapsto z$ and $e_{1}\mapsto e_{1}+z$. Using the fact that $z$ is central in $\hcV\jump{z}$, we then observe that 
\begin{align}
c_{z}(h)
&=c_{z}(e_{\infty}+f)\\
&=e_{\infty}-z+c_{z}(f)\\
&=e_{\infty}-z+f\\
&=h-z,
\end{align}
where we used $f\in\fG$ in the third equality.
From this computation and the assumption, we have
\begin{align}
(h-z)c_z(Q)^{-1}
&=c_z(hQ^{-1})\\
&=c_{z}(Pe_{\infty})\\
&=c_{z}(P)(e_{\infty}-z).
\end{align}
Then we consider the images under the operator
\[\rho_{z\to e_{\infty}}\colon\hcV\jump{z}\simeq\hcV\hotimes R\jump{z}\stackrel{z\mapsto e_{\infty}}{\longrightarrow}\hcV\hotimes\hcV\stackrel{\text{mult.}}{\longrightarrow}\hcV.\]
It is easily checked that $\rho_{z\to e_{\infty}}(c_{z}(P)(e_{\infty}-z))=0$. Thus we obtain the equality $h\psi_{0}=\psi_{0}e_{\infty}$ where $\psi_{0}\coloneqq\rho_{z\to e_{\infty}}(c_{z}(Q)^{-1})$. We also see that this series $\psi_{0}$ is invertible in $\hcV$. Hence the element
\[\psi \coloneqq \psi_{0}\cdot \psi_{0}(0,e_{\infty})^{-1},\]
satisfies the equation $\psi e_{\infty}\psi^{-1}=h$. To obtain the group-likeness of $\psi$, we compute the element
\[d \coloneqq (\psi\otimes\psi)^{-1}\Delta^{\cV}(\psi)\in\hcV.\]
From the assumption on $h$, we have
\begin{align}
\psi e_{\infty}\psi^{-1}\otimes 1+1\otimes\psi e_{\infty}\psi^{-1}
&=h\otimes 1+1\otimes h\\
&=\Delta^{\cV}(h)\\
&=\Delta^{\cV}(\psi e_{\infty}\psi^{-1})\\
&=(\psi\otimes\psi)d(e_{\infty}\otimes 1+1\otimes e_{\infty})d^{-1}(\psi\otimes\psi)^{-1},
\end{align}
so that $d$ commutes with $\Delta^{\cV}(e_{\infty})=e_{\infty}\otimes 1+1\otimes e_{\infty}$. Since the centralizer of $\Delta^{\cV}(e_{\infty})$ in $\hcV^{\hotimes 2}$ is equal to $R\jump{e_{\infty}}^{\hotimes 2}$, we conclude that $d=1$.
\end{proof}

\begin{lemma}\label{lem:gpt_fact}
Let $a,b,c\in\hcV$. Then the following are equivalent:
\begin{enumerate}
\item $c=be_{\infty}a+a+b$.
\item The condition (1) of Lemma \ref{lem:schneps_fact} holds for
\[h = e_{\infty} + e_{0}ae_{0} + e_{1}ce_{0} + e_{1}be_{1}.\]
\end{enumerate}
\end{lemma}
\begin{proof}
The implication (1) $\implies$ (2) is easy, as there is an explicit expression $(P,Q)=(1-e_1b,1-ae_0)$. 
We show the rest. Decomposing the given equation as
\begin{align}
h &= Pe_{\infty}Q\\
&= (1+e_0P^0+e_1P^1)e_{\infty}(1+Q_{0}e_0+Q_1e_1)\\
&= \begin{multlined}[t]e_{\infty} + e_{0}(P^0e_{\infty}Q_0-P^0-Q_0)e_{0} + e_{1}(P^1e_{\infty}Q_0 - P^1 - Q_0)e_{0} \\+ e_1(P^1e_{\infty}Q_1-P^1-Q_1)e_1 + e_{0}(P^0e_{\infty}Q_{1}-P^0-Q_1)e_1,\end{multlined}
\end{align}
we obtain
\[(a,b,c) = (P^0e_{\infty}Q_0-P^0-Q_0,P^1e_{\infty}Q_1-P^1-Q_1,P^1e_{\infty}Q_0 - P^1 - Q_0),\]
and
\[P^0e_{\infty}Q_1-P^0-Q_1 = 0.\]
Hence we have $Q_{1}=(P^0e_{\infty}-1)^{-1}P^0$ and thus
\begin{align}
    a &= (P^{0}e_{\infty}-1)Q_{0} - P^{0} \\
    &= (P^{0}e_{\infty}-1)(Q_{0} - Q_{1}), \\
    b &= (P^1e_{\infty}-1)Q_1-P^{1} \\
    &= (P^1e_{\infty}-1)(P^0e_{\infty}-1)^{-1}P^{0}-P^{1}\\
    &= (P^1e_{\infty}-1)P^{0}(e_{\infty}P^{0}-1)^{-1}-P^{1}\\
    &= (P^{1}-P^{0})(e_{\infty}P^0-1)^{-1}.
\end{align}
Hence we obtain
\begin{align}
be_{\infty}a + a + b
&= (P^{1}-P^{0})(e_{\infty}P^0-1)^{-1}e_{\infty}(P^{0}e_{\infty}-1)(Q_{0} - Q_{1})+a+b\\
&= (P^{1}-P^{0})e_{\infty}(Q_{0} - Q_{1})+a+b \\
&= P^{1}e_{\infty}Q_{0} + P^{0}e_{\infty}Q_{1} - P^0-Q_0 - P^1-Q_1 \\
&= c.
\end{align}
It completes the proof.
\end{proof}

\begin{lemma}\label{lem:red_fact}
Assume that $\varphi\in\cG$. Then the following are equivalent:
\begin{enumerate}
\item $\varphi\in\cG^{\mu}_{\Red}$.
\item The condition (1) of Lemma \ref{lem:schneps_fact} holds for
\[h = - (1+e_1B_{\mu}(e_1))\varphi e_0 - e_1\varphi(1+B_{\mu}(e_0)e_{0}).\]
\end{enumerate}
\end{lemma}
\begin{proof}
From \eqref{eq:reduction_condition}, we know that (1) is equivalent to
\begin{align}E_{\mu}(\varphi) &\coloneqq \varphi(B_{\mu}(e_{0})-\varphi^{0})+(B_{\mu}(e_{1})-\varphi_{1})\varphi+\varphi_{0}(1+e_{0}\varphi^{0})+(1+\varphi_{1}e_{1})\varphi^{1}\\ &=0.\end{align}
Now we put
\begin{align}
a &= -\varphi^0,\\
b &= -\varphi_1,\\
c &= -\varphi_0-\varphi^1-B_{\mu}(e_1)\varphi-\varphi B_{\mu}(e_0),
\end{align}
so that $h=e_{\infty}+e_{0}ae_{0}+e_{1}ce_{0}+e_{1}be_{1}$.
Then we see that the equality
\begin{align}
E_{\mu}(\varphi)
&= -c + (\varphi_0e_0-\varphi)\varphi^{0} + \varphi_{1}(e_1\varphi^1-\varphi)\\
&= -c - (1+\varphi_1e_1)\varphi^{0} - \varphi_{1}(1+e_0\varphi^0)\\
&= -c + be_{\infty}a + a + b
\end{align}
holds. Hence $\varphi\in\cG^{\mu}_{\Red}$ if and only if the equality $c=be_{\infty}a+a+b$ holds. Therefore, from Lemma \ref{lem:gpt_fact}, we obtain the desired equivalence.
\end{proof}

\begin{lemma}\label{lem:inert_fact}
Assume that $\varphi\in\cG$. Then the following are equivalent:
\begin{enumerate}
\item $\varphi\in\cG^{\mu}_{\inert}$.
\item The condition (1) of Lemma \ref{lem:schneps_fact} holds for
\[h\coloneqq -\bch_{\mu}(e_1,\varphi e_0\varphi^{-1}).\]
\end{enumerate}
\end{lemma}
\begin{proof}
    First we rephrase (1); let $\varphi\in\cG^{\mu}_{\inert}$ so that we have an element $\psi\in\exp(\hff)$ satisfying
    \[ \bch_{\mu}(\varphi e_{0}\varphi^{-1},\psi e_{\infty}\psi^{-1}) = -e_{1}. \]
    Then, using \cite[Lemma 14.7]{ef25}, we have
    \begin{align}
        h   
        &= -\bch_{\mu}(e_{1},\varphi e_{0}\varphi^{-1},\psi e_{\infty}\psi^{-1},-\psi e_{\infty}\psi^{-1})\\
        &= -\bch_{\mu}(e_{1},-e_{1},-\psi e_{\infty}\psi^{-1})\\
        &= \psi e_{\infty}\psi^{-1}.
    \end{align}
     Hence we see that the condition $\varphi\in\cG^{\mu}_{\inert}$ is equivalent to $h=\psi e_{\infty}\psi^{-1}$ for some $\psi\in\exp(\hff)$. Since Lemma \ref{lem:schneps_fact} holds, it is also equivalent to that $h=Pe_{\infty}Q$ for some $P$ and $Q$ and their constant terms are $1$.
\end{proof}

\begin{theorem}\label{thm:main_first_honbun}
    We have the following.
    \begin{enumerate}
        \item $\fG_{\Red}=\fG_{\inert}$.
        \item $\cG^{\mu}_{\Red}=\cG^{\mu}_{\inert}$ for all $\mu\in R$.
    \end{enumerate}
\end{theorem}
\begin{proof}
    The assertion (1) is same as Lemma \ref{lem:ef16}. For (2), it suffices to show the equivalence between Lemma \ref{lem:schneps_fact} (1) for 
    \[h \coloneqq -(1+e_1B_{\mu}(e_1))\varphi e_0 - e_1\varphi(1+e_0B_{\mu}(e_0))\]
    and the same one for
    \[h' \coloneqq -\bch_{\mu}(e_1,\varphi e_0\varphi^{-1}).\]
    We consider the formal power series
    \[\ep_{\mu}(t) \coloneqq \frac{\exp(\mu t)-1}{\mu t} \in R\jump{t}.\]
    Since the equation
    \[(1+t B_{\mu}(t))\ep_{\mu}(t) = \frac{\exp(\mu t)+1}{2}\]
    holds by definition, we have
    \begin{align}
        &\ep_{\mu}(e_{1})h\ep_{\mu}(e_{0})\\
        &=-\ep_{\mu}(e_{1})((1+e_{1}B_{\mu}(e_{1}))\varphi e_{0}+e_{1}\varphi (1+e_{0}B_{\mu}(e_{0})))\ep_{\mu}(e_{0}) \\
        &= -\frac{e^{\mu e_{1}}+1}{2}\varphi e_{0}\ep_{\mu}(e_{0}) - \ep_{\mu}(e_{1})e_{1}\varphi\frac{e^{\mu e_{0}}+1}{2} \\
        &= -\frac{1}{2\mu}\left((e^{\mu e_{1}}+1)\varphi(e^{\mu e_{0}}-1) + (e^{\mu e_{1}}-1)\varphi (e^{\mu e_{0}}+1)\right) \\
        &= -\frac{1}{\mu}\left(e^{\mu e_{1}}\varphi e^{\mu e_{0}}-\varphi\right) \\
        &= -\frac{1}{\mu}(\exp(-\mu h')-1)\varphi\\
        &= h'\ep_{\mu}(-h')\varphi.
    \end{align}
    Note that this computation is always available even if $\mu$ is not invertible in $R$. By this identity
    \[\ep_{\mu}(e_{1})h\ep_{\mu}(e_{0}) = h'\ep_{\mu}(-h')\varphi,\]
    the condition (1) of Lemma \ref{lem:schneps_fact} for $h$ is equivalent to that of $h'$, because all of the factors
    \[\ep_{\mu}(e_{1}), \ep_{\mu}(e_{0}),\ep_{\mu}(-h'),\varphi\]
    are invertible.
\end{proof}

\part{Flexion theory and Kashiwara--Vergne relations}
In this part, we exhibit an interpretation of the Kashiwara--Vergne structure as a explicit family among multiple zeta values. Our tool used here is Ecalle's flexion theory. For the basic theory, see \cite{ecalle11}, \cite{schneps15}, \cite{fhk23}, \cite{kawamura25} and \cite{kawamura25}. We basically adopt with conventions in \cite{kawamura25}. Especially, a bimould is often regarded as a linear function on the algebra $\fV$ of biwords ($=$ double-layered sequence of formal variables) equipped with the Hopf algebra structure via the shuffle product and deconcatenation coproduct.
\section{Dimorphic transportation}
\subsection{Decks and families of functions}
    First we remark that, some nice examples of \emph{families of functions} \cite[\S1]{fhk23} are recovered as a \emph{deck} \cite[\S2]{kawamura251}. Here we show an example used later: consider the deck $\cX_{\ser}$ over $\bbQ$ defined by
    \[r \mapsto \cX_{\ser,r}\coloneqq(R\mapsto \cF_{\ser}(V_{r}\otimes R)),\]
    where the family $\cF_{\ser}$ is the one in \cite[Example 3]{fhk23}. Similarly, we define the deck $\cX_{\Lau}$ with respect to the family $\cF_{\Lau}$ in \cite{fhk23}. Though its image is not representable as an endofunctor of $\CAlg_{\bbQ}$, but we can check that it does not prevent any argument done in this paper.

\subsection{Translation map}
For a given deck $\cX_{\bullet}$ coming from a named family of functions, set
\[
    \BIMU^{\bu\text{-}\bullet,\vcst}\subseteq\BIMU(R;\cX_{\bullet})
\]
to be, as in \cite{kawamura261}, the space of all bimoulds independent of variables of the lower layer. We have the isomorphism of algebras
    \begin{align}
        \ma\colon\cV&\to\BIMU^{\user,\vcst};~\varphi \mapsto\ma(\varphi),\\
        \ma(\varphi)\binom{u_{1},\ldots,u_{r}}{v_{1},\ldots,v_{r}}
        &=\sum_{k_{1},\ldots,k_{r}\ge 1}
        \langle\varphi\mid (-e_{1})e_{0}^{k_{1}-1}\cdots
        (-e_{1})e_{0}^{k_{r}-1}\rangle\nonumber\\
        &\qquad\cdot u_{1}^{k_{1}-1}\cdots
        (u_{1}+\cdots+u_{r})^{k_{r}-1}.
    \end{align}
    which is defined in \cite[\S 3]{schneps15}. It restricts to the isomorphism of completed Lie algebras
    \[\ma \colon (\hff,\langle,\rangle) \to \ARI^{\user,\vcst}.\]
    Moreover, we define the subspace $\ARI^{\ucom,\vcst}\subset \ARI^{\user,\vcst}$ by adding the vanishing condition of degree $1$:
    \[\ARI^{\ucom,\vcst} \coloneqq \{A\in\ARI^{\user,\vcst}\mid A\binom{0}{0}=0\},\]
    as well as subspaces or subalgebras of $\ARI^{\user,\vcst}$.
\subsection{Double shuffle and Kashiwara--Vergne structures in flexion theory}
    Let us denote by $\dmr_{0}$ (resp.~$\krv_{2}$) the corresponding Lie algebra to the pro-unipotent group $\DMR_{0}$ (resp.~$\KRV_{2}$). They are Lie subalgebras of $(\fG,\langle,\rangle)$.
    The following are shown in \cite{schneps15}, \cite{fkrs}:
    \begin{theorem}[{\cite{fkrs}, \cite{schneps15}}]
        By restricting $\ma$, we get the following isomorphisms.
        \begin{enumerate}
            \item $\dmr_{0} \simeq \ARI_{\baral\ast\baril}^{\ucom,\vcst}$.
            \item $\krv_{2} \simeq \ARI_{\al+\senu\ast\pusinu}^{\ucom,\vcst}$.
        \end{enumerate}
    \end{theorem}
    Moreover, regarding the depth-graded variants, we have the following diagram drawing the so-called dimorphic transportations:
    \begin{theorem}[{\cite{kawamura25}}]
        This diagram of Lie algebras commutes:
        \[\xymatrix{
            \ARI_{\baral\ast\baral}^{\uLau,\vcst} \ar[r]^-{\adari(\pal)}_-{\simeq} & \ARI_{\baral\ast\barfol}^{\uLau,\vcst} \\ 
            \ARI_{\al+\push\ast\pusnu}^{\uLau,\vcst} \ar[r]^-{\adari(\pal)}_-{\simeq} \ar@{^{(}->}[u]^-{\text{incl.}} & \ARI_{\al+\Esena\ast\Opusnu}^{\uLau,\vcst} \ar@{^{(}->}[u]_-{\text{incl.}}
        }\]
    \end{theorem}

\section{Weak alternality and alternility} 
\subsection{Quasi-shuffle product with flexion units}
We use the product $\sh_{\ast}$ in a generalized sense. It is based on Komiyama's idea\footnote{Though it originates in \cite{komiyama21}, he kindly shared his more detailed note about it during Kagoshima Mould Seminar 2025.}. Let $\fO$ be a flexion unit. Consider the $\bbQ$-linear tensor product $\fV_{\cO}\coloneqq \fV\otimes_{\bbQ}\cX_{1}(\bbQ)$ of the biword space $\fV$ with the valued space $\cX_{1}(\bbQ)$ of $\cO$. We define the quasi-shuffle product $\sh_{\ast}$ associated with $\cO$ by the recursive rule
\begin{multline}
    \alpha\bw\sh_{\ast}\beta\bw' = \alpha(\bw\sh_{\ast}\beta\bw') + \beta(\alpha\bw\sh_{\ast}\bw') - \cO(\alpha\flr{\beta})(\bw\sh_{\ast}\ful{\alpha}\beta\bw') - \cO(\fll{\alpha}\beta)(\alpha\fur{\beta}\bw\sh_{\ast}\bw'),
\end{multline}
which is valid for arbitrary biwords $\bw,\bw'$ and ones $\alpha,\beta$ of length $1$. It is $\fV_{\cO}$-bilinear and has $1$ (the empty word) as its unit.
\subsection{Definition and translation}
\begin{definition}
    Let $A\in\LU$ be a bimould. Then we say that $A$ is \emph{weakly alternal} (resp.~\emph{weakly alternil}) if
    \[A(\bw\sh(w'_{1}))=0\quad (\text{resp. }A(\bw \sh_{\ast} (w'_{1}))=0),\]
    for every non-empty biword $\bw\in\fV$ and $w'_{1}$ (biword of length $1$) where the unit is taken as $\fO=\Pi$.
\end{definition}
As the alternality (resp.~alternility) is usually abbreviated as $\al$ (resp.~$\il$) in the notation, we denote by $\al^{\weak}$ (resp.~$\il^{\weak}$) for the weak alternality (resp.~weak alternility).

\subsection{Transportation of weak dimorphy}
\begin{lemma}
    The Lie algebra automorphism $\adari(\pal)\in\Aut(\ARI^{\uLau,\vcst})$ restricts to the isomorphism
    \[\adari(\pal)\colon\ARI_{\baral\ast\baral^{\weak}}^{\uLau,\vcst} \stackrel{\simeq}{\longrightarrow} \ARI_{\baral\ast\baril^{\weak}}^{\uLau,\vcst}.\]
\end{lemma}
\begin{proof}
    From Lemma \ref{lem:main_fourth_flexion}, we see that $\preari(-,B)=\preira(-,B)$ for every $B\in\ARI_{\al\ast\al^{\weak}}$. Hence we can show
    \[\adari(\pal)(\ARI_{\baral\ast\baral^{\weak}}^{\uLau,\vcst})\subseteq\ARI_{\baral\ast\baril^{\weak}}^{\uLau,\vcst},\]
    in the same way as \cite[Theorem 4.6.1]{schneps15}. The converse argument is almost same as \cite[Theorem 3.5]{kawamura25}.
\end{proof}

\subsection{Weak double shuffle relations}
Recall that the torsor $\DMR_{\mu}$ is defined as the subset of all $\varphi\in\cG$ such that $\langle\varphi\mid e_{0}e_{1}\rangle=\mu^{2}/24$ and $\varphi_{\Ra}^{\ast}\in\cW$ is primitive with respect to the harmonic coproduct $\Delta_{\ast}$. The latter condition is rewritten as the so-called harmonic relation
\[\zeta_{\varphi}^{\ast}(w \ast w') = \zeta_{\varphi}^{\ast}(w)\zeta_{\varphi}^{\ast}(w'),\]
for all words $w,w'$ in the alphabet $\{y_{k}\}_{k\ge 1}$.
\begin{definition}
    \begin{enumerate}
        \item For $\mu\in R$, we define the \emph{weak double shuffle set} $\WDMR_{\mu}$ as the subset of $\cG$ consisting of $\varphi\in\cG$ such that $\langle\varphi\mid e_{0}e_{1}\rangle=\mu^{2}/24$ and the condition
        \[\zeta_{\varphi}^{\ast}(w \ast y_{l}) = \zeta_{\varphi}^{\ast}(w)\zeta_{\varphi}^{\ast}(y_{l}),\]
        satisfied for every word $w$ in the alphabet $\{y_{k}\}_{k\ge 1}$ and every positive integer $l$.
        \item We denote by $\wdmr_{0}$ the corresponding Lie algebra to the subgroup $\WDMR_{0}\subseteq\DMR_{0}$.
    \end{enumerate}
\end{definition}
Recall that two coproduct structures, $\Sh$ and $\Sh_{\ast}$, of $\BIMU^{\user,\vcst}$ are defined in \cite[\S 2.1]{komiyama21}. The bimould $A\in\LU^{\user,\vcst}$ is alternal (resp.~alternil) if and only if
\[A\otimes 1+1\otimes A=\Sh(A)\quad (\text{resp. }=\Sh_{\ast}(A)),\]
as shown in \cite[Propositions 2.4, 2.7]{komiyama21}.
\begin{lemma}\label{lem:wdmr_flexion}
    The map $\ma$ restricts to the Lie algebra isomorphism
    \[\ma\colon \wdmr_{0} \stackrel{\simeq}{\longrightarrow} \ARI_{\baral\ast\baril^{\weak}}^{\ucom,\vcst}.\]
\end{lemma}
\begin{proof}
    It immediately follows from the fact that the translation map $\mathsf{mi}\coloneqq\swap\circ\ma$ intertwines $\Sh_{\ast}$ and $\Delta^{\cW}$ \cite[Lemma 131]{fhk23}.
\end{proof}
\begin{lemma}\label{lem:krv_flexion}
     The map $\ma$ restricts to the Lie algebra isomorphism
    \[\ma\colon \krv_{2} \stackrel{\simeq}{\longrightarrow} \ARI_{\baral\ast\baril^{\weak}}^{\ucom,\vcst}.\]
\end{lemma}
\begin{proof}
    Same as \cite[Corollary 4.20]{fkrs}.
\end{proof}
\subsection{Proof of Theorem \ref{thm:main_fourth}}
\begin{lemma}\label{lem:explicit_shuffle}
    Let $\bw$ be a biword in $\fV$ and $w_{0}$ a biword of length $1$. Then we have
    \[A((w_{0})\sh \bw) = \sum_{i=0}^{r}(\push^{-1}\circ\pus^{-i}\circ\push^{i+1})(A)(w_{0}\bw)\]
    for any bimould $A$.
\end{lemma}
\begin{proof}
    Denote $P_{i} \coloneqq \push^{-i}\circ\pus^{i-1}\circ\push$. Then we have
    \[P_{i}(A)(a_{1}\cdots a_{r}) = A(a_{i}a_{1}\cdots a_{i-1}a_{i+1}\cdots a_{r})\]
    for every biword $a_{1}\cdots a_{r}\in \fV$. Hence we obtain
    \begin{align}
        A((w_{0})\sh \bw)
        &=\sum_{i=0}^{r}A(w_{1}\cdots w_{i}w_{0}w_{i+1}\cdots w_{r})\\
        &=\sum_{i=0}^{r}P_{i+1}^{-1}(A)(w_{0}w_{1}\cdots w_{r}),
    \end{align}
    which is the desired equation.
\end{proof}
\begin{lemma}\label{lem:explicit_push_invariance}
    For a biword $\bw=w_{1}\cdots w_{r}$ ($r\ge 2$) and an alternal bimould $A$, we have
    \begin{align}
        (\push-\id)(A)(\bw)
        &=(-1)^{r-1}\swap(A)\left(
        \binom{-v_{r}}{-u_{r}}\sh
        \binom{v_{1},v_{2}-v_{1},\ldots,v_{r-1}-v_{r-2}}
              {u_{1}+\cdots+u_{r-1},\ldots,u_{r-1}}
        \right)\nonumber\\
        &\quad-\swap(A)\left(
        \binom{v_{r}}{u_{1}+\cdots+u_{r}}\sh
        \binom{v_{r-1}-v_{r},\ldots,v_{1}-v_{2}}
              {u_{1}+\cdots+u_{r-1},\ldots,u_{1}}
        \right).
    \end{align}
\end{lemma}
\begin{proof}
    Put
    \begin{align}
        a &= \binom{v_{r}}{u_{1}+\cdots+u_{r}},\\
        b_{i} &= \binom{v_{r-i}-v_{r+1-i}}{u_{1}+\cdots+u_{r-i}},\\
        c &= \binom{-v_{r}}{-u_{r}},\\
        d_{i} &= \binom{v_{i}-v_{i-1}}{u_{i}+\cdots+u_{r-1}}
    \end{align}
    for $1\le i\le r-1$. Here we regard $v_{0}$ as $0$. Since $\swap(A)(ab_{1}\cdots b_{r-1})=A(\bw)$, we have
    \[\swap(A)(a \sh b_{1}\cdots b_{r-1}) = A(\bw) + \sum_{i=1}^{r-1}\swap(A)(b_{1}\cdots b_{i}ab_{i+1}\cdots b_{r-1}).\]
    On the other hand, since the equalities
    \[(\swap\circ\anti)(A)(cd_{1}\cdots d_{r-1}) = \push(A)(\bw)\]
    and
    \[(\swap\circ\anti)(A)(d_{1}\cdots d_{r-i}cd_{r-i+1}\cdots d_{r-1}) = \swap(A)(b_{1}\cdots b_{i}ab_{i+1}\cdots b_{r-1})\]
    hold for $1\le i\le r-1$ by definition, we obtain
    \begin{align}
        (\swap\circ\anti)(A)(c \sh d_{1}\cdots d_{r-1}) &= \push(A)(\bw) + \sum_{i=1}^{r-1}\swap(A)(b_{1}\cdots b_{i}ab_{i+1}\cdots b_{r-1})\\
        &=\push(A)(\bw) -A(\bw)+\swap(A)(a \sh b_{1}\cdots b_{r-1}).
    \end{align}
    Since we assumed that $A$ is alternal, it holds that
    \[(\swap\circ\anti)(A)(c \sh d_{1}\cdots d_{r-1}) = (-1)^{r-1}\swap(A)(c \sh d_{1}\cdots d_{r-1}),\]
    due to the $\mantar$-invariance \cite[Lemma 2.5.3]{schneps15}. 
\end{proof}
\begin{remark}
    Lemma \ref{lem:explicit_push_invariance} also gives a refinement of the argument done in \cite[Lemma 3.7--Proposition 3.11]{fk}.
\end{remark}
\begin{lemma}\label{lem:main_fourth_flexion}
    We have the equality
    \[\ARI_{\al+\push\ast\pusnu} = \ARI_{\baral\ast\baral^{\weak}}.\]
\end{lemma}
\begin{proof}
    The $\subseteq$ part follows from Lemma \ref{lem:explicit_shuffle}. For the converse inclusion, the $\push$-invariance follows from Lemma \ref{lem:explicit_push_invariance}, and then the $\pus$-neutrality can be shown by using Lemma \ref{lem:explicit_shuffle} again because $\push$ is bijective.
\end{proof}
\begin{lemma}\label{lem:wdmr_group}
    The functor $R\mapsto\WDMR_{0}$ forms a pro-unipotent $\bbQ$-group with the product $\circledast$.
\end{lemma}
\begin{proof}
    It follows from the fact that $(\DMR_{0},\circledast)$ forms a group and that $\Delta^{\cW}$ and $\circledast$ is compatible with the depth filtration of $\hcW$.
\end{proof}
\begin{lemma}\label{lem:wdmr}
    We have the equality
    \[\exp_{\circledast}(\wdmr_{0}) = \WDMR_{0}.\]
\end{lemma}
\begin{proof}
    Since $\WDMR_{0}$ is a subgroup of $(\cG,\circledast)$ and $\wdmr_{0}$ is obtained by the linearization of the group functor $R\mapsto\WDMR_{0}$, we have the result.
\end{proof}
\begin{lemma}\label{lem:intersection_wdmr_solkv}
    We denote by $M_{\mu}$ the set of all associators in $\hcV$ with its coupling constant $\mu\in R$. Then
    \begin{enumerate}
        \item $M_{\mu}\subseteq \WDMR_{\mu}$.
        \item The injection $\nu^{\EF}$ restricts to $M_{\mu}\hookrightarrow \SolKV_{\mu}$.
    \end{enumerate}
\end{lemma}
\begin{proof}
    The second assertion is shown in \cite[\S 7]{at12}. For the first, combine the result $M_{\mu}\subseteq \DMR_{\mu}$ \cite[Corollary 0.4]{furusho11} with the definition of $\WDMR_{\mu}$.
\end{proof}
\begin{theorem}[$=$ Theorem \ref{thm:main_fourth}]\label{thm:main_fourth_honbun}
    We have the following.
    \begin{enumerate}
        \item $\krv_{2} = \wdmr_{0}$.
        \item The restriction of $\nu^{\EF}$ gives $\WDMR_{\mu} \stackrel{\simeq}{\longrightarrow}\SolKV_{\mu}$.
    \end{enumerate}
\end{theorem}
\begin{proof}
    The first assertion is a consequence of Lemmas \ref{lem:main_fourth_flexion}, \ref{lem:krv_flexion} and \ref{lem:wdmr_flexion}. We shall show the second.
    First we check that
    \[\KRV_{2} = \exp(\krv_{2}) = \nu^{\EF}(\exp_{\circledast}(\wdmr_{0})) = \nu^{\EF}(\WDMR_{0}),\]
    where the first (resp.~third) equality is a consequence of Proposition \ref{prop:racinet_composition} (resp.~Lemma \ref{lem:wdmr})
    We also know that there is a non-empty subset $S$ of $\nu^{\EF}(\WDMR_{\mu})\cap\SolKV_{\mu}$ due to Lemma \ref{lem:intersection_wdmr_solkv}. Since both of $(\KRV_{2},\SolKV_{\mu})$ and $(\WDMR_{0},\WDMR_{\mu})$ forms a torsor, we have
    \begin{align}
        \nu^{\EF}(\WDMR_{\mu})
        &=\nu^{\EF}(\WDMR_{0}\circledast (\nu^{\EF})^{-1}(S))\\
        &=\nu^{\EF}(\WDMR_{0})\cdot S\\
        &=\KRV_{2}\cdot S\\
        &=\SolKV_{\mu}.
    \end{align}
    This completes the proof.
\end{proof}

\begin{corollary}\label{cor:wdmr_reduction}
    The pair $(\WDMR_{0},\WDMR_{\mu})$ is a subtorsor of $(\cG^{0}_{\Red},\cG^{\mu}_{\Red})$.
\end{corollary}

\part{Stability for the $S_{3}$-symmetry}
\section{Three involutions}
Enriquez--Furusho \cite{ef25} also defined an important involution $\Theta^{\mu}$ acting on the torsor $\cG^{\mu}_{\inert}$. It originates from the involutive automorphism of the curve $\bbP^{1}\setminus\{0,1,\infty\}$ intertwining $0$ and $\infty$. Here we shall construct the other two, namely transpositions of $\{1,\infty\}$ and $\{0,1\}$.
\begin{lemma}[{\cite[Lemma 0.17 (a)]{ef25}}]\label{lem:partner}
    For any element $\varphi\in\cG^{\mu}_{\inert}$, we find the unique element $\psi\in\exp(\hff_{2})$ satisfying
    \begin{enumerate}
        \item $\bch_{\mu}(\varphi e_{0}\varphi^{-1},\psi e_{\infty}\psi^{-1})=-e_{1}$.
        \item $\log(\psi)\equiv \mu e_{1}/2$ mod $\overline{[\hff,\hff]}$.
    \end{enumerate}
\end{lemma}
For a given $\varphi\in\cG^{\mu}_{\inert}$, the element $\psi$ indicated in Lemma \ref{lem:partner} is called the \emph{partner} of $\varphi$. The partner is also denoted as $h_{\varphi}$.

\begin{definition}\label{def:involutions}
    Suppose that $\varphi\in\cG^{\mu}_{\inert}$ and let $\psi$ be its partner. We define
    \begin{align}
        \Theta^{\mu}(\varphi) & \coloneqq e^{-\mu e_{1}/2}\psi(e_{\infty},e_{1}),\\
        \rT^{\mu}(\varphi) & \coloneqq e^{-\mu e_{1}/2}\psi(e_{0},e_{\infty})^{-1}e^{\mu e_{\infty}}\varphi(e_{0},e_{\infty})e^{-\mu e_{0}/2},
    \end{align}
    and
    \[\Delta^{\mu}(\varphi) \coloneqq \varphi(e_{1},e_{0})^{-1}.\]
\end{definition}
\begin{remark}
    Obviously the definition of $\Delta^{\mu}$ is available on the whole $\cG$, not necessarily on $\cG^{\mu}_{\inert}$, and we do not need write $\mu$ explicitly in the notation $\Delta^{\mu}$.
\end{remark}
\begin{lemma}\label{lem:delta_racinet}
    The map $\Delta^{\mu}$ preserves the Racinet product.
\end{lemma}
\begin{proof}
    Taking two elements $\varphi_{1}$ and $\varphi_{2}\in\cG$ arbitraily, we see that
    \begin{align}
        &\Delta^{\mu}(\varphi_{1}\circledast\varphi_{2})\\
        &=\Delta^{\mu}(\varphi_{2}(\varphi_{1}e_{0}\varphi_{1}^{-1},e_{1})\varphi_{1})\\
        &=\varphi_{1}(e_{1},e_{0})^{-1}\cdot\varphi_{2}(\varphi_{1}(e_{1},e_{0})e_{1}\varphi_{1}(e_{1},e_{0})^{-1},e_{0})^{-1}\\
        &=\varphi_{2}(e_{1},\varphi_{1}(e_{1},e_{0})^{-1}e_{0}\varphi_{1}(e_{1},e_{0}))^{-1}\varphi_{1}(e_{1},e_{0})^{-1},
    \end{align}
    holds, which is equal to $\Delta^{\mu}(\varphi_{1})\circledast\Delta^{\mu}(\varphi_{2})$.
\end{proof}
\begin{lemma}\label{lem:delta_inert}
    The map $\Delta^{\mu}$ induces an involutive automorphism of the torsor $(\cG^{0}_{\inert},\cG^{\mu}_{\inert})$.
\end{lemma}
\begin{proof}
    Since it is clearly involutive and Lemma \ref{lem:delta_racinet} holds, we show $\cG^{\mu}_{\inert}$ is stable under the action of $\Delta^{\mu}$. Take an element $\varphi\in\cG^{\mu}_{\inert}$ and let $\psi$ be its partner. Then we have
    \begin{align}
        &\bch_{\mu}(\Delta^{\mu}(\varphi)e_{0}\Delta^{\mu}(\varphi)^{-1},\Delta^{\mu}(\varphi)e^{-\mu e_{0}}\Delta^{\mu}(\psi)^{-1}e_{\infty}(\Delta^{\mu}(\varphi)e^{-\mu e_{0}}\Delta^{\mu}(\psi)^{-1})^{-1})\\
        &=\bch_{\mu}(\varphi(e_{1},e_{0})^{-1}e_{0}\varphi(e_{1},e_{0}),\varphi(e_{1},e_{0})^{-1}e^{-\mu e_{0}}\psi(e_{1},e_{0})e_{\infty}\psi(e_{1},e_{0})^{-1}e^{\mu e_{0}}\varphi(e_{1},e_{0}))\\
        &=\varphi(e_{1},e_{0})^{-1}e^{-\mu e_{0}}\cdot \bch_{\mu}(e_{0},\psi(e_{1},e_{0})e_{\infty}\psi(e_{1},e_{0})^{-1})\cdot e^{\mu e_{0}}\varphi(e_{1},e_{0}).
    \end{align}
    Applying the transposition $e_{0}\leftrightarrow e_{1}$ on the formula
    \[\bch_{\mu}(e_{1},\psi e_{\infty}\psi^{-1}) = -e^{\mu e_{1}}\varphi e_{0}\varphi^{-1}e^{-\mu e_{1}},\]
    which is easily obtained by definition of $\psi$, we obtain
    \begin{align}
        &\bch_{\mu}(\Delta^{\mu}(\varphi)e_{0}\Delta^{\mu}(\varphi)^{-1},\Delta^{\mu}(\varphi)e^{-\mu e_{0}}\Delta^{\mu}(\psi)^{-1}e_{\infty}(\Delta^{\mu}(\varphi)e^{-\mu e_{0}}\Delta^{\mu}(\psi)^{-1})^{-1})\\
        &=-\varphi(e_{1},e_{0})^{-1}e^{-\mu e_{0}}\cdot e^{\mu e_{0}}\varphi(e_{1},e_{0})e_{1}\varphi(e_{1},e_{0})^{-1}e^{-\mu e_{0}}\cdot e^{\mu e_{0}}\varphi(e_{1},e_{0})\\
        &=-e_{1}.
    \end{align}
    Hence we see that $\Delta^{\mu}(\varphi)\in\cG^{\mu}_{\inert}$.
\end{proof}
\begin{corollary}\label{cor:s3}
    The subset
    \[\{\id,\Theta^{\mu},\rT^{\mu},\Delta^{\mu},\Delta^{\mu}\circ\Theta^{\mu},\Delta^{\mu}\circ\rT^{\mu}\}\]
    forms a subgroup of $\Aut(\cG^{\mu}_{\inert})$ (automorphism group of the $\cG^{0}_{\inert}$-torsor) which is isomorphic to $S_{3}$.
\end{corollary}
\begin{proof}
    As we already showed that $\Theta^{\mu}$ and $\Delta^{\mu}$ are involutive, it suffices to show the equality
    \[\rT^{\mu} = \Delta^{\mu}\circ\Theta^{\mu}\circ\Delta^{\mu}.\]
    The equality
    \[\bch_{\mu}(\Delta^{\mu}(\varphi)e_{0}\Delta^{\mu}(\varphi)^{-1},\Delta^{\mu}(\varphi)e^{-\mu e_{0}}\Delta^{\mu}(\psi)^{-1}e_{\infty}(\Delta^{\mu}(\varphi)e^{-\mu e_{0}}\Delta^{\mu}(\psi)^{-1})^{-1})=-e_{1},\]
    obtained in the proof of Lemma \ref{lem:delta_inert}, we see that
    \begin{align}
        h_{\Delta^{\mu}(\varphi)}
        &=\Delta^{\mu}(\varphi)e^{-\mu e_{0}}\Delta^{\mu}(\psi)^{-1}e^{-\mu e_{\infty}/2}\\
        &=\varphi(e_{1},e_{0})^{-1}e^{-\mu e_{0}}\psi(e_{1},e_{0})e^{-\mu e_{0}/2}.
    \end{align}
    It yields
    \[(\Theta^{\mu}\circ\Delta^{\mu})(\varphi) = e^{-\mu e_{1}/2}\varphi(e_{1},e_{\infty})^{-1}e^{-\mu e_{\infty}}\psi(e_{1},e_{\infty})e^{-\mu e_{0}/2}.\]
    Taking the inverse of the right-hand side and then applying the transposition $e_{0}\leftrightarrow e_{1}$, we get
    \begin{align}
        (\Delta^{\mu}\circ\Theta^{\mu}\circ\Delta^{\mu})(\varphi)
        &= e^{-\mu e_{1}/2}\psi(e_{0},e_{\infty})^{-1}e^{\mu e_{\infty}}\varphi(e_{0},e_{\infty})e^{-\mu e_{0}/2}\\
        &=\rT^{\mu}(\varphi).
    \end{align}
\end{proof}

\section{Ren's result revisited}
\subsection{$\Theta$-invariance and multitangent functions}
Ren \cite{ren26} obtained the result claming the equivalence between the $\Theta^0$-invariance and the so-called \emph{harmonic antipode} relation, in the Lie algebra setting. In the torsor setting, the harmonic antipode is desrcibed as follows:
\begin{lemma}[{\cite[p.~58]{hoffman}}]
  The tuple $(\hcW,\cdot,\Delta^{\cW})$ forms a completed Hopf $R$-algebra with the antipode defined by
  \[S^{\ast}\colon \varphi \mapsto \rev_{\ry}(\varphi(e_{0}+e_{1},-e_{1})).\]
\end{lemma}
Ren's statement is easily generalized to the torsor setting, namely the following theorem is obtained in the same way as \cite[Theorem 1.2]{ren26}:
\begin{theorem}\label{thm:ren_enhanced}
  For any element $\mu\in R^{\times}$ and $\varphi\in\cG_{\inert}^{\mu}$, the following are equivalent:
  \begin{enumerate}
    \item $S^{\ast}(\varphi_{\Ra}^{\ast})\varphi_{\Ra}^{\ast} = 1$.
    \item $\Theta^{\mu}(\varphi)=\varphi$.
  \end{enumerate}
\end{theorem}
\begin{proof}
  Quite same as \cite{ren26}.
\end{proof}
Since $\DMR_{\mu}$ is characterized as the intersections of group-like parts of $\hcV$ and $\hcW$, the following corollary is immediate.
\begin{corollary}\label{cor:dmr_theta}
  The set $\DMR_{\mu}$ has the inclusion into the $\Theta^{\mu}$-invariant part of $\cG_{\inert}^{\mu}$.
\end{corollary}
We now give a brief reinterpretation of this corollary: assume $\mu=0$ for simplicity. Let $\varphi\in\DMR_{0}$. Then, combining Theorems \ref{thm:main_fourth_honbun} and \ref{thm:main_first_honbun}, we obtain
\[\varphi\in\DMR_{0}\subseteq\WDMR_{0}=\KRV_{2}\subset\cG^{0}_{\inert}=\cG^{0}_{\Red}.\]
Hence the equality $\cT_{\varphi}=\cR_{\varphi}$ holds. Since the series $\cT_{\varphi}$ is group-like in the completed Hopf algebra $\hcW\hotimes R(\!(z)\!)$, the antipode formula $S^{\ast}(\cT_{\varphi})=\cT_{\varphi}^{-1}$ holds. In particular, if $r\ge 1$ the equality
\[\sum_{i=0}^{r} (-1)^{r-i}\hzeta(k_{1},\ldots,k_{i};z)\hzeta^{\star}(k_{r},\ldots,k_{i+1};z)=0\]
for any $k_{1},\ldots,k_{r}\ge 1$, where $\hzeta^{\star}$ is given by
\[\hzeta^{\star}(l_{1},\ldots,l_{s};z) \coloneqq \sum_{\substack{c\ge 1\\ j_{1},\ldots,j_{c}\ge 1\\ j_{1}+\cdots+j_{c}=s}}\hzeta(l_{1}+\cdots+l_{j_{1}},\ldots,l_{j_{1}+\cdots+j_{c-1}+1}+\cdots+l_{s};z).\]
Comparing the coefficients of $z^{-2}$ and rewriting it in terms of $\varphi$, we obtain
\[ e_{1} + \varphi e_{0}\varphi^{-1} + \varphi(e_{\infty},e_{1})e_{\infty}\varphi(e_{\infty},e_{1})^{-1} = 0, \]
which asserts $\Theta^{0}(\varphi)=\varphi$.
\subsection{Towards duality deduction problem}
Among the research community of multiple zeta values, the following problem, the so-called \emph{duality deduction problem}, is long-standing:
\begin{quote}
    Would regularized double shuffle relations imply the duality relation?
\end{quote}
Since the duality relation for MZVs is equivalent to $\Delta^{2\pi i}(\Phi_{\KZ})=\Phi_{\KZ}$, the problem is rewritten as
    \begin{equation}\label{eq:ddp}
    \Delta^{\mu}|_{\DMR_{\mu}} \stackrel{?}{=} \id
    \end{equation}
in terms of the involution $\Delta^{\mu}$. Here we show a new criterion for this problem. Since the equality $\SolKV_{\mu}\stackrel{?}{=}\DMR_{\mu}$ is expected, we can propose a much weaker version:
\begin{conjecture}\label{conj:ddp}
    The weak double shuffle relations would imply the harmonic antipode relation. Namely,
    \[\varphi\in\WDMR_{\mu} \stackrel{?}{\Longrightarrow} S^{\ast}(\varphi_{\Ra}^{\ast})=(\varphi_{\Ra}^{\ast})^{-1}.\]
\end{conjecture}
\begin{theorem}[$=$ Theorem \ref{thm:main_extra}]\label{thm:main_extra_honbun}
    If Conjecture \ref{conj:ddp} is true, the weak double shuffle relations imply the duality relation. In particular, under this assumption, the duality deduction problem \eqref{eq:ddp} is true.
\end{theorem}
\begin{proof}
    Assume that Conjecture \ref{conj:ddp} is true and take $\varphi\in\WDMR_{\mu}$. 
    Since $\SolKV_{\mu}$ is shown in \cite{at12} to be stable under the action of $S_{3}$ (Corollary \ref{cor:s3}), the torsor $\WDMR_{\mu}$ is stable too, via Theorem \ref{thm:main_fourth_honbun}. Therefore we have $\rT^{\mu}(\varphi)\in\WDMR_{\mu}$. From the assumption, we obtain
    \[S^{\ast}((\rT^{\mu}(\varphi))_{\Ra}^{\ast})\cdot (\rT^{\mu}(\varphi))_{\Ra}^{\ast}=1,\]
    and hence $(\Theta^{\mu}\circ\rT^{\mu})(\varphi)=\rT^{\mu}(\varphi)$ by Theorem \ref{thm:ren_enhanced}. 
    It yields that
    \[\Delta^{\mu}(\varphi) = (\rT^{\mu}\circ\Theta^{\mu}\circ\rT^{\mu})(\varphi)=\varphi,\]
    which means $\Delta^{\mu}|_{\WDMR_{\mu}}=\id$.
\end{proof}

\section*{Acknowledgments}
The authors thank Professor Hidekazu Furusho and Professor Shin-ichiro Seki for their valuable comments and discussions. They are also grateful to Dr.~Takumi Anzawa for reading the manuscript carefully.

\section*{Use of AI}
The authors used ChatGPT 6.2 Astra to help refine the proof of Theorem~\ref{thm:main_first_honbun} and identify relevant references.
The authors verified the contents of the paper and take full responsibility for the final manuscript.

\end{document}